\documentclass[12pt,leqno]{article} 

\usepackage{amsmath} 
\usepackage{amscd} 
\usepackage{amssymb} 
\usepackage{amsfonts}
\usepackage{amsthm}
\usepackage{verbatim} 
\usepackage{pictex}
\rm 

\newlength{\fixboxwidth}
\newtheorem{thm}{Theorem}
\newtheorem{cor}{Corollary}
\newtheorem{rem}{Remark}
\newtheorem{definition}{Definition}
\newtheorem{lemma}{Lemma}
 
\newtheorem{proposition}{Proposition}

\def\wt{\widetilde }

\def\d{\delta }
\def\e{\varepsilon } 
\def\epsilon{\varepsilon } 
 
\def\rho{\varrho } 
 
\def\phi{\varphi }
\def\a{\alpha }
\def\b{\beta }
\def\({\biggl( } 
\def\){\biggr) }

\def\sq2{{\sqrt{2}}} 
\def\sgn{{\rm sgn }} 
 
\def\A{{\mathcal A}} 
  
\def\lin{{\rm lin}}
\def\non{{\rm non}} 
\def\cont{{\rm cont}} 
\def\cF{\mathcal{F}} 
\def\frame{{\rm frame}}

\def\cB{\mathcal{B}}
\def\cF{\mathcal{F}}
\def\cG{\mathcal{G}}
\def\cK{\mathcal{K}}
\def\cD{\mathcal{D}}
\def\cs{\mathcal{S}}

\newcount\refnum
\def\refx{\smallskip \global\advance\refnum by 1 {[\the\refnum ] \ }}

\def\b{\beta}
\def\a{\alpha}
\def\d{\delta}

\def\L{{\cal L}} 
\def\N{{\cal N}} 

\def\ce{{\mathcal E}}
\def\cl{{\mathcal L}}
\def\cx{{\mathcal X}}
\def\C{{\mathcal C}}

\def\D{{\cal D}}
\def\S{{\cal S}} 
\def\Nb{{\mathbb N}} 
\def\R{{\mathbb R}}

\def\tor{{\mathbb T}} 
\def\Rd{{\mathbb R}^d} 
 
\def\Z{{\mathbb Z}} 
\def\Zd{{\mathbb Z}^d}

\def\O{{\mathcal O}}

\def\supp{{\rm supp \, }}

\def\dist{{\rm dist \, }}

\def\lsim{\raisebox{-1ex}{$~\stackrel{\textstyle <}{\sim}~$}}

\newcommand{\eproof}{\qquad \hfill  \qedsymbol \\}
 
\title{Optimal Approximation of Elliptic Problems by Linear 
and Nonlinear Mappings III: Frames}

\author{Stephan Dahlke\thanks{This author acknowledges support through 
the European Union's Human Potential Programme, under contract
HPRN-CT-2002-00285 (HASSIP), and through DFG, 
Grant Da 360/4-3. He also wants to thank the
Friedrich-Schiller-Universit\"at Jena for the 
hospitality and support.}, Erich Novak, Winfried Sickel \\ 
\hbox{\small dahlke@mathematik.uni-marburg.de, novak@math.uni-jena.de,
sickel@math.uni-jena.de} }

\begin{document}

\maketitle

\begin{center} 
{\large \em Dedicated to our dear colleague and friend Henryk Wo\'zniakowski 
on the occasion of his 60th birthday.}
\end{center} 

\begin{abstract}
We study the optimal approximation of the solution of an operator 
equation
$ 
\A(u) = f 
$ 
by certain $n$-term approximations with 
respect to specific classes of frames. 
We consider 
worst case errors, where $f$ is an element of the unit ball of a 
Sobolev or Besov space
$B^t_q(L_p(\Omega))$ and $\Omega \subset \R^d$ is a bounded
Lipschitz domain; the error is always measured in the $H^s$-norm. 
We study the order of convergence of the 
corresponding nonlinear frame widths
and compare it with  several other approximation schemes. 
Our main result is that the approximation order is the same as for    
the nonlinear widths associated with  Riesz bases, the Gelfand widths,
and the manifold widths. 
This order is better than the order of the 
linear widths iff $p<2$.
The main advantage of frames compared to Riesz bases, 
which were studied in our earlier papers, is the fact that we
can now handle arbitrary bounded Lipschitz domains---also 
for the upper bounds. 

\end{abstract}




\section{Introduction}         


We study the optimal approximation of the solution of an operator 
equation
\begin{equation}     \label{e01} 
\A(u) = f , 
\end{equation} 
where $\A$ is a linear operator 
\begin{equation}     \label{e02} 
\A: H  \to G
\end{equation} 
from a Hilbert space $H$ to another Hilbert space $G$. 
We always assume that $\A$ is boundedly invertible, hence 
\eqref{e01} has a unique solution for any $f \in G$.  
We have in mind the more specific situation of an 
operator equation which is given as follows. 
Assume that $\Omega \subset \R^d$ is a bounded Lipschitz 
domain and assume that 
\begin{equation}     \label{e03} 
\A : H^s_0 (\Omega) \to H^{-s} (\Omega) 
\end{equation} 
is an isomorphism, where $s > 0$. 
For the exact definitions of Lipschitz domains and spaces 
of distributions defined on such domains we refer to the Appendix, see also \cite{dns2}.
Now we put $H=H^s_0(\Omega)$ and $G = H^{-s} (\Omega)$. 
Since $\A$ is boundedly invertible, the inverse mapping 
$S: G \to H$ is well defined. This mapping is sometimes 
called the solution operator---in particular if we want 
to compute the solution $u=S(f)$ from the given 
right-hand side $\A(u)=f$. 

We  study different  mappings $S_n$   for the 
approximation of the solution $u=\A^{-1}(f)$  for $f$ contained in    
$F \subset G$. We consider the worst case error 
\begin{equation}   \label{e04} 
e(S_n, F,H) = 
\sup_{\Vert f \Vert_F \le 1} \Vert \A^{-1}(f)- S_n(f) \Vert_H , 
\end{equation}  
where $F$ is a normed (or quasi-normed) space, $F \subset G$. 
In our main results, $F$ will be a Sobolev or Besov 
space.\footnote{Formally we only deal with Besov spaces. 
Because of the embeddings $B_1^{-s+t}(L_p(\Omega)) \subset 
W_p^{-s+t}(\Omega)  \subset B_\infty^{-s+t}(L_p(\Omega))$, 
which hold for $1 \le p \le \infty$, $t\ge s$, see \cite{T02}, 
our results are valid also 
for Sobolev spaces.} 
Hence we use the following 
commutative diagram
\begin{eqnarray*}
G & \stackrel{\hbox{$S$}}{\longrightarrow} & H \\
I \,  & \nwarrow \qquad  \nearrow & S_F \\
& F. &
\end{eqnarray*}
Here $I: \, F \to G$ denotes the identity and $S_F$ the restriction of 
$S$ to $F$.
Then one is interested in approximations that have an optimal order of
convergence  depending on $n$, 
where $n$ denotes the degrees of freedom.
For our purposes,  the following approximation schemes are important. 
Consider the class $\L_n$ of all continuous linear mappings 
$S_n : F \to H$, 
$$ 
S_n(f) = \sum_{i=1}^n L_i(f) \cdot \tilde  h_i 
$$ 
with arbitrary $\tilde h_i \in H$.  
The worst case error of optimal linear 
mappings is given by the \emph{approximation numbers} 
or \emph{linear widths}
$$ 
e_n^\lin (S , F , H) = \inf_{S_n \in \L_n}  e(S_n, F,H) . 
$$ 
We may also use nonlinear approximations with respect to a
Riesz basis $\mathcal R$  of $H$, i.e., 
we consider the class $\N_n (\mathcal R)$ of all (linear or 
nonlinear) mappings of the form 
$$  
S_n(f) = \sum_{k=1}^n c_k \,  h_{i_k} ,      
$$ 
where the $c_k$ and the $i_k$ depend in an arbitrary way on $f$. 
Then the \emph{nonlinear widths}  
$e_{n,C}^\non (S,F,H)$ are given by  
$$ 
e_{n,C}^\non (S, F , H) = \inf_{\mathcal R \in {\mathcal R}_C} \inf_{S_n \in \N_n(\mathcal R)} 
e(S_n, F,H). 
$$ 
Here ${\mathcal R}_C$ denotes a set of Riesz bases for $H$  
where $C$ indicates the stability of the basis, i.e.,
we require $B/A\le C$
and $A,B$ are the Riesz constants of the basis.
The investigation of these widths $e_{n,C}^\non$ and 
its comparison with the linear widths have been the major 
part of our analysis in \cite{dns1,dns2}.
This has continued earlier research on related topics, cf. e.g.
\cite{Ka85,Tem00,Tem02,Tem03}.
The next type of widths we are interested 
in has served as a very useful tool
in our analysis of the widths $e_{n,C}^\non$ in \cite{dns2}.
The \emph{manifold widths} are related
to the class  $\C_n$  of continuous mappings, given by arbitrary 
continuous mappings
$N_n : F \to \R^n$ and $\phi_n : \R^n \to H$. 
Again we define the worst case error of optimal continuous 
mappings by
\begin{equation}\label{manifold}
e_n^{\rm cont} (S, F , H) = \inf_{S_n \in \C_n}  e(S_n, F, H ) , 
\end{equation}
where $S_n = \phi_n \circ N_n$.  
These numbers  have been studied in 
\cite{DHM89,Ma90} and later in 
\cite{dns2,DKLT93,DD96,DD00}.
As mentioned above we have studied the relationships 
of these widths in \cite{dns2}.
It has turned out that 
for problems as in  \eqref{e03} with 
$F=B^{-s+t}_q(L_p(\Omega))$ 
(with some extra conditions on $\Omega$) 
one has the following:
if $p \ge 2$ and $t>0$ then
\begin{eqnarray} 
e_n^{\rm lin} (S, B^{-s+t}_q(L_p(\Omega)) , H^s_0 (\Omega))
\asymp &
\!\!\! \!\!\!\!\!\! 
\!\!\!\!\!\! \!\!\!\!\!\! 
e_n^{\rm cont} (S, B^{-s+t}_q(L_p(\Omega)) , H^s_0 (\Omega))
\\ 
\asymp & 
e_{n,C}^{\rm non} (S, B^{-s+t}_q(L_p(\Omega)) , H^s_0 (\Omega)) 
\asymp n^{-t/d}\, ,  \nonumber 
\end{eqnarray} 
whereas in the case $0 < p <2$ with $t>d(1/p - 1/2)$
\[
e_n^{\rm lin} (S, B^{-s+t}_q(L_p(\Omega)) , H^s_0 (\Omega))
\asymp n^{-t/d + 1/p-1/2}
\]
and
\[ 
e_n^{\rm cont} (S, B^{-s+t}_q(L_p(\Omega)) , H^s_0 (\Omega))
\asymp
e_{n,C}^{\rm non} (S, B^{-s+t}_q(L_p(\Omega)) , 
H^s_0 (\Omega)) \asymp n^{-t/d}\, .
\]
Hence, if $p<2$ then there is an essential difference 
in the behavior, nonlinear approximations can 
do better than linear ones. 

This paper is a continuation 
of \cite{dns1,dns2}. We are again
interested in  optimal nonlinear  approximation 
schemes, but this time not
related to  Riesz bases but to classes of frames. 
The motivation for this
is given by the following observations.  
In \cite{dns2}, we presented upper
and lower bounds for $e_{n,C}^\non (S, F , H)$.  
The proof of the lower bound
was quite general and used the fact that  
$e_{n,C}^\non (S, F , H)$ can be
estimated  from below by the manifold 
widths $e_n^{\rm cont} (S, F , H)$ up to some 
constants. In contrary to this, the
proof of the upper bound was based on  
norm equivalences of Besov norms with
weighted sequence norms that are induced by a biorthogonal wavelet
basis. However, this restricts the 
choice of the underlying domain $\Omega
\subset \R^d$ since on a general Lipschitz 
domain the construction of a suitable wavelet
basis might be very complicated or even impossible. 
This problem becomes less
serious in the frame setting since a suitable 
wavelet frame always  exists, see
Section \ref{pbasic1} for a detailed discussion. 
Moreover, in recent years the application of
frame methods for the numerical resolution of the solution $u$ in 
(\ref{e01})
has become  a field of increasing importance. 
Especially, it has been possible
to derive adaptive wavelet frame schemes 
that are guaranteed to converge for a
wide range of problems \cite{DFR04,DFRSW06,stev03}. 
Therefore it is important to
clarify the power that frame schemes can  have, in principle. 

In this paper, we give a first answer. Our main result states that the
nonlinear frame widths show the same asymptotic behavior as the 
$ e_{n,C}^\non (S, F , H)$, where we now can allow arbitrary 
bounded Lipschitz domains. 

There is an interesting difference 
to the  Riesz bases case.  In the frame setting, 
we do \emph{not} work with arbitrary $n$-term approximations, 
but only with those induced 
by a frame pair, see Section
\ref{introfw} for details.  
The reason is that, for practical applications, only these 
canonical representations are used. 
Actually we prove that if we would allow arbitrary 
$n$-term approximations then 
the  associated frame widths would be
zero. Moreover, certain
conditions related to stability must be satisfied by the 
admissible frames. Fortunately,
these conditions are always satisfied for the 
known constructions of wavelet frames on Lipschitz domains. 

This  paper is organized as follows.  
In Section \ref{setting}, we describe the
basic setting. First of all,  we 
introduce and discuss
the frame concept as far as it is needed for our 
purposes. Then, in Subsection
\ref{introfw}, we define the 
nonlinear frame widths and prove some
basic properties that are needed in the sequel.  Section
\ref{jetzt} contains the main results 
of this paper. 
In the next section two examples are discussed: the Poisson equation for 
Lipschitz domains and  a Fredholm integral equation of the first kind
(the single layer potential).
Proofs of our main results  
are given in Section \ref{Beweis}. For general
Hilbert spaces $H$ and $G$ we show that similar 
to the  Riesz bases case the
nonlinear frame widths can be estimated from below by the  manifold
widths. Then, for the  more specific case of Besov spaces on Lipschitz
domains, we  also prove  an upper estimate 
which shows that the asymptotic
behavior is the same as for the Riesz basis case---but this
time for arbitrary bounded Lipschitz domains. 

{\bf Notation.} 
We write $a \asymp b$ if there exists  a constant $c>0$ 
(independent of the context dependent 
relevant parameters) such that
\[
c^{-1} \, a \le b \le c \, a \,.
\]
One-sided estimates of this type are denoted 
by $a \lsim b$. 
All unimportant constants will be denoted by $c$, 
sometimes with additional indices. 
Identity operators are always denoted by $I$, 
also sometimes with additional indices.


\section{Frames} \label{setting}


In this paper, we  will study certain 
approximations of $u=S(f)$ based on
frames. Therefore, in this section we  
recall the  basic properties of frames
as far as they are needed for our purposes and introduce  the corresponding nonlinear widths.  
For further information on frames,
we refer the reader  e.g. to \cite{Chr03, Gro00}.    
A sequence $\mathcal{F}=\{h_k\}_{k \in \Nb}$ in  a separable Hilbert
space ${H}$ is a \textit{frame} for ${H}$ if 
there exist constants $A, B >0$ such that 
\begin{equation}
\label{framestability}
A^2 \sum_{k=1}^\infty \bigl|( f, h_k)_{H}\bigr|^2
\le \|f\|_{H}^2 \leq  
B^2 \sum_{k=1}^\infty \bigl|( f, h_k)_{H}\bigr|^2
\end{equation}
for all $f\in {H}$.
As a consequence of \eqref{framestability}, the corresponding 
operators of analysis and synthesis
given by
\begin{equation}
\label{analysisop}
{\mathcal T}:{H}\rightarrow\ell_2(\Nb),
\quad f\mapsto\bigl(( f,h_k)_{H}\bigr)_{k\in\Nb},
\end{equation}
\begin{equation}
\label{synthesisop}
{\cal T}^*:\ell_2(\Nb)\rightarrow {H},
\quad\mathbf{c}\mapsto\sum_{k = 1}^\infty c_kh_k,
\end{equation}
are bounded.
The composition ${\mathcal T}^*{\mathcal T}$ is a 
boundedly invertible (positive and self-adjoint) operator 
called the {\it frame operator}. Furthermore,  
$\widetilde{\mathcal{F}}
:= ({\mathcal T}^*{\mathcal T})^{-1}\mathcal{F}$ is 
again a frame for ${H}$, the {\it canonical dual frame}.
The following formulas hold
\begin{equation}
\label{reconstruction}
f  = \sum_{k =1}^\infty 
( f, ({\mathcal T}^*{\mathcal T})^{-1} h_k)_H \, 
h_k = \sum_{k =1}^\infty ( f, h_k )_H \, 
({\mathcal T}^*{\mathcal T})^{-1} h_k
\end{equation}
for all $f\in H$.
This classical concept of a frame is 
too general, we need an additional stability  condition, 
stronger than \eqref{framestability}.
Without this additional assumption on the frames,
there 
would not exist lower bounds for corresponding widths  as we shall now explain. 

\begin{rem}   \label{rem1} 
Let $H$ be a separable Hilbert space and let $K \subset H$ be a 
compact subset. Then for an arbitrary $C>1$ 
there exists a frame 
$\mathcal{F}=\{h_i\}_{i \in \Nb}$ in 
$H$ with $B/A <C$ such that the following is true: 
For all $f \in K$ and for all $\e >0$ there exists a 
$h_i \in \mathcal{F}$ and $c \in \R$ such that 
$$
\Vert f - c h_i \Vert_H <\e .
$$
Hence the best $n$-term approximation yields an error 0 
already for $n=1$. 
To prove this statement, we construct such a frame for a given
compact set $K \subset H$. 
Let $M_1= \{ e_i, \ i \in \Nb \}$ be a complete orthonormal 
set of $H$ and let $\{ k_i , \ i \in \Nb \}$ be a dense subset
of $K$. 
We consider sets of the form 
$$
M_2^\delta  = \{ \a_1 k_1, \, \a_2 k_2 , \, \dots \} \subset H
$$
with $\a_i = \delta^i$, where $0< \delta <1$  and put 
$\mathcal{F_\delta} = M_1 \cup M_2^\delta$. 
It is not difficult to check that 
$\mathcal{F_\delta}$ is a frame with all the claimed 
properties if $\d = \d (C)$ is chosen appropriately. 

The frames 
$\mathcal{F_\delta}$ can be considered as ``pathological'',
since the norms of many elements of 
$\mathcal{F_\delta}$ are extremely small. A first idea would 
be to request that  the norms of the frame elements are uniformly bounded from above and below, 
\[
0<c_1 \leq \Vert h_i \Vert_H \leq c_2<\infty \qquad \mbox{for all}\quad 
h_i \in \mathcal{F}=\{h_i\}_{i \in \Nb}\, ,
\]
but this does not
help: 
Now we can define $\mathcal{F_\delta}$ as the union of 
$M_1$ and multiples of the $e_i \pm \a_i k_i$. 
Then one obtains 
such a ``normed'' frame such that:  
For all $f \in K$ and for all $\e >0$ there exist
$h_i \in \mathcal{F}$ and $c_i \in \R$ such that 
$$
\Vert f - c_1 h_1 - c_2 h_2  \Vert_H <\e .
$$
Therefore we go into a different direction, see Definitions~{\rm \ref{defpair} }
and {\rm \ref{defiwidth}}. 
\end{rem} 


\subsection{Frame Pairs} \label{albtraum}


As it is well-known, Sobolev spaces built on $L_2 (\Omega)$ 
can be discretized by means of weighted $\ell_2$-spaces, 
see the Appendix for some examples how one can do this.
Let $w:= (w_k)_{k\in \Nb}$ be a sequence of positive 
numbers which we call simply a {\em weight} in what follows. Then we put
\[
\ell_{2,w}:= \Big\{a= (a_k)_{k \in \Nb}: \quad 
\| \, a \, \|_{\ell_{2,w}} := \Big(\sum_{k=1}^\infty
w_k \, |a_k|^2\Big)^{1/2}<\infty 
\Big\}\, .
\]

\begin{definition} \label{defpair}  
Let $H$ be  a separable Hilbert space with dual space $H'$. 
Let $w= (w_k)_k$ be a weight.
\\
{\rm (i)}
Two sequences  $({\mathcal F}, {\mathcal G})$, 
${\mathcal F}:=\{h_k\}_{k \in {\Nb}} \subset H', 
{\mathcal G}:=\{g_k\}_{k \in {\Nb}} \subset H$, 
are  called a  {\em frame pair} for $(H,w)$, if
\begin{equation} \label{atom}
f=\sum_{k=1}^ {\infty}{\langle f,h_k  \rangle}_{H\times H'}\, g_k, 
\end{equation}
holds for all $f \in H$ and   we have the norm equivalence 
\begin{equation} \label{normpair}
A \,  \|\, (\langle f, h_k\rangle_{H\times H'})_{k \in \Nb}\, 
\|_{\ell_{2,w}}  \le 
 \|f\|_H \le 
B \,  \|\, (\langle f, h_k\rangle_{H\times H'})_{k \in \Nb}\, 
\|_{\ell_{2,w}} 
\end{equation}
with some positive constants $A, B$. In addition, we require that
there exists a bounded linear operator $R~:~\ell_{2,w}\longrightarrow H$  satisfying
\begin{equation} \label{blablabla}
R(\delta_k)=g_k\quad \mbox{and} \quad  \| R\| \leq B.
\end{equation}
\\
{\rm (ii)}
Let  $K$ be a subspace of $H$.  
A frame pair $({\mathcal F}, {\mathcal G})$ for $(H,w)$ 
is called {\em stable} with respect to $K$ if the inequality
\begin{equation} \label{star}
A'\,  \|(
\langle f, h_k\rangle_{H\times H'})_{k\in \Lambda} \|_{\ell_{2,w}}
\le 
\left\|\sum_{k \in \Lambda} \langle f, h_k\rangle_{H\times H'} \, 
g_k \right\|_H 
\end{equation}
holds with some $A'>0$, all finite subsets  $\Lambda \subset {\Nb}$ and all
$f \in K$. \\ 
{\rm (iii)} Let  $K$ be a  subspace  of $H$ and 
let $C\ge 1$ be a given number.
By  ${\cal P}_C(K)$  we
denote the set of all stable frame pairs  $({\cal G}, {\cal F})$  with 
respect to $K$  
such that the constants $A, B$ and $A'$ in {\rm (\ref{normpair})} and {\rm (\ref{star})} 
satisfy $B /\min(A,A') \le C$.
\end{definition}

\begin{rem} \label{notation}
To avoid any type of confusion we shall 
use $(\cdot, \cdot)$ for the scalar product in 
$H$ and $\langle \cdot , \cdot\rangle$ 
for duality pairings, in particular for 
$H \times H'$.
\end{rem}

Some comments are in order. 

\begin{rem} \label{ganzviellaber}
\begin{itemize}
\item[(i)]  A frame pair in the sense of {\rm (\ref{atom})} and {\rm (\ref{normpair})} is sometimes called
an {\rm atomic decomposition}, cf. e.g. {\rm \cite[Def.~17.3.1.]{Chr03}}. However,
the phrase {\em atomic decomposition} is used with a different meaning in
the theory of function spaces, cf. e.g. {\rm \cite{FJ90,Ky03,T92,T06a}}.
For this reason we do not use it here.
\item[{(ii)}] Let $({\mathcal F}, 
{\mathcal G})$ be a frame pair for $(H,w)$.
As above let
${\mathcal F}=\{h_k\}_{k \in {\Nb}} 
\subset H'$ and ${\mathcal G}:=\{g_k\}_{k \in {\Nb}} \subset H$.
By the Riesz representation theorem, for 
every $h_k$ there exists an element
$\widetilde{h}_k \in H$ such that  
$\langle f, h_k \rangle_{H\times H'}=(f,
\widetilde{h}_k)_H$. 
Consequently,
\[
\|\, ( f, \, \sqrt{w_k}\, 
\widetilde{h}_k)_{k \in \Nb}\, \|_{\ell_{2}}
=
\|\, (\langle f, h_k\rangle_{H\times H'})_{k \in \Nb}\, 
\|_{\ell_{2,w}} \qquad \mbox{for all} \quad  f \in H\, .
\]
Hence, there is a one-to-one correspondence 
between $\cF$ and the Hilbert frame 
$(\sqrt{w_k}\, \widetilde{h}_k)_k$. 
However, note that $\cG$ need not be related to the 
canonical dual frame of $(\sqrt{w_k}\, \widetilde{h}_k)_k$.
\item[(iii)] The reader might wonder why we use the concept of frame pairs instead of the classical frame
setting as introduced in {\rm (\ref{framestability})} and {\rm (\ref{reconstruction})}.  However, since  we
 are  dealing here with Gelfand triples 
$(H^s_0(\Omega), L_2(\Omega), H^{-s}(\Omega))$, $s-1/2 \neq \mbox{integer}$, see
Remark \ref{tilde}, this
approach would be at least problematic since we are not allowed to identify the space
$H^s_0(\Omega)$ with its dual.(Otherwise, it would not be possible to identify $L_2(\Omega)$
 with its dual at the same time - a strange construction. We refer to {\rm \cite{hack}} for further details.)
\item[(iv)] Our concept is closely related to 
Banach frames in the sense of
{\rm  \cite{Gro91,Gro04}}. A \textit{Banach frame} for a separable and 
reflexive Banach space $\cB$ is a sequence 
$\cF =\{h_k\}_{k \in \Nb}$ in $\cB'$ with 
an associated sequence space 
$\cB_d$ such that the following properties hold:
\setlength{\leftmargin}{4cm}
\begin{itemize}
\item[{\rm (B1)}] norm equivalence: 
there exist constants $A,B>0$  such that 
\begin{equation}
\label{bfnormequiv}
A \,  \left\|\bigl(\langle f,h_k
\rangle_{\cB\times\cB'}\bigr)_{k
\in \Nb}\right\|_{\cB_d} 
\le \|f\|_{\cB} \le 
B \,  \left\|\bigl(\langle f,h_k
\rangle_{\cB\times\cB'}\bigr)_{k
\in \Nb}\right\|_{\cB_d} 
\end{equation} 
for all $f \in \cB$;
\item[{\rm (B2)}] there exists a bounded 
operator $\mathcal S $ from $\cB_d$ onto $\cB$, a so-called 
{\em synthesis} or {\em reconstruction operator}, such that
\begin{equation}
\label{banachframeop}
{\mathcal S} \left(\bigl(\langle f,h_k\rangle_{\cB\times\cB'}
\bigr)_{k \in \Nb}\right) =f.
\end{equation}
\end{itemize}
(It is a remarkable fact that for Banach spaces the existence of the
reconstruction operator does not follow from the norm equivalence
{\rm (\ref{bfnormequiv})}  
and has to be explicitly required). 

A  frame pair in the sense of Definition {\rm  \ref{defpair} (i)}  
induces a Banach
frame ${\cal F}=\{h_k\}_{k\in \Nb}$ for the special  case $\cB=H$,
$\cB_d=\ell_{2,w}(\Nb)$ where the operator $R$  serves as synthesis operator,  cf. {\rm \cite[Thm.~3.2.3]{Chr03}}.
Consequently, in our setting, 
the estimate
\begin{equation} \label{upperwichtig}
\left\|\sum_{k \in \Nb}c_k \, g_k\right\|_{H} 
\le B \, \left\| (c_k)_{k \in \Nb}\right\|_{\ell_{2,w}}
\end{equation}
always holds.  
\item[{(v)}]
We comment on the condition  {\rm  (\ref{star})}.
Clearly, {\rm (\ref{star})} always holds on all of $H$  
for a Riesz basis  $\{g_k\}_{k \in \Nb} $ for $H$. 
However, there exist frames which are not Riesz bases and for which 
{\rm (\ref{star})} holds on $H$. E.g.,
take an orthonormal basis $\{e_k\}_{k \in \Nb}$ and define the frame
$\cF=:\{e_1, 2^{-1/2}e_2, 2^{-1/2}e_2, e_3, e_4...\}$.
This is a tight frame, \eqref{normpair} holds 
with $A=B=1$, so the primal and the canonical 
dual frame coincide. (We
refer again to {\rm \cite[Chapt.~5]{Chr03}} for further information). 
Since  $\{e_k\}_{k \in \Nb}$ is an orthonormal basis, 
a direct computation shows that  {\rm (\ref{star})} 
holds for $A'=2^{-1/2}$.
Nevertheless,  requiring {\rm (\ref{star})} on all of $H$ 
would be very restrictive, and most frames would not satisfy it. As an
example, consider the frame $\cF:=\{e_1, 2^{-1/2}e_2, 2^{-1/2}e_2,
3^{-1/2}e_3, 3^{-1/2}e_3,  3^{-1/2}e_3,\ldots \}$. 
This is also  a tight
frame, but again a direct check shows 
that {\rm (\ref{star})}
does not hold. Therefore we require {\rm (\ref{star})} 
only on  subsets.
Fortunately, such a  condition is satisfied in case 
of the known frame constructions
for function spaces on Lipschitz domains.

\item[(vi)]  The example in (v)  shows that the two constants $A$ amd $A'$ in
Definition {\rm  \ref{defpair}} need not to be related at all.  Nevertheless, to
avoid unnecessary notational difficulties, we will restrict ourselves to the
case $A=A'$ in the sequel. The modifications to the case $A\not= A'$ are
straightforward. 
\item[(vii)]  For simplicity, we have introduced our basic concepts for frame pairs  indexed by the set of
natural numbers.  Later on, we shall also use frame pairs corresponding to more general countable sets,
with the obvious modifications.
\end{itemize}
\end{rem}
For later use, let us finally state the following 
simple but useful property: frame 
pairs are invariant under isomorphic mappings.

\begin{lemma} \label{framemap}  
Let $G,H$  be Hilbert  spaces and 
let 
$S : G \to H$ be  an isomorphism. 
Let $(\cF, {\mathcal G})$  be a frame pair  for $(G,w)$
with frame constants  $A,B$. Then the following holds:\\
{\rm (i)} $({{S}^*}^{-1}(\cF), {S}({\mathcal G}))$  is a  frame
pair   for $(H,w)$ with frame constants
$\wt A = A/\|S^{-1} \|$ and
$\wt B= B\|{S}\|$.  \\
{\rm (ii)}  If  $(\cF, {\mathcal G})$ 
is contained in ${\mathcal P}_C(K)$
then  $({{S}^*}^{-1}(\cF), S({\mathcal G}))$ is contained in
${\mathcal P}_{\wt C}(S(K))$, where 
$\wt C=C\|{S}\|\|{S}^ {-1}\|$.
\end{lemma}

\begin{proof} {\em Step 1.} Proof of (i).  
We start by showing (\ref{atom}).  For $f \in H$, we obtain
\begin{eqnarray*} f=S(S^{-1}(f))=S 
\Big( \sum_{k\in \Nb} \langle S^{-1}(f),
h_k\rangle_{H\times H'}\,  g_k \Big) = \sum_{k \in \Nb} \langle f, 
{S^*}^{-1}(h_k)\rangle_{H\times H'} \, S(g_k).
\end{eqnarray*}
The next step is to  show the norm equivalence  (\ref{normpair}). We obtain
\begin {eqnarray*} \frac{1}{\|S\|}\|f\|_H &=&  
\frac{1}{\|S\|}\|(S\circ S^{-1})(f)\|_H\leq \|{S}^{-1}(f)\|_G \\
&\le & B  (\langle S^{-1}(f),
 h_n\rangle_{G\times G'})_{\ell_{2,w}} = 
B  (\langle f, {{S}^*}^{-1}(h_k)\rangle_{H\times H'})_{\ell_{2,w}}
\\
&\le & \frac B A \, \|S^{-1}(f)\|_G\leq \frac B A\,  \|S^{-1}\| \|f\|_H \, .
\end{eqnarray*}
Let $R$ be the bounded operator associated with $(\cF, {\mathcal G})$. Then 
$\tilde{R}=S \circ R$ is again a bounded operator with
$$\tilde{R}(\delta_k)=S(R(\delta_k))=S(g_k), \quad \|\tilde{R}\|\leq \|S\| \|R\|\leq \|S\| B,$$
and (i) is shown. \\
{\em Step 2.} Proof of (ii).
For $f \in S(K)$, we get
\begin{eqnarray*}
\|\sum_{k \in \Lambda} \langle f, {{S}^*}^{-1}(h_k)
\rangle_{H\times H'} \, S(g_k)\|_{H} &\geq& \|S^{-1}\|^{-1}\Big\| 
\sum_{k \in \Lambda} \langle
S^{-1}(f), h_k\rangle_{G \times G'} \, g_k\, \Big\|_G
\\
& \geq & \|S^{-1}\|^{-1}A\, \|\, 
(\langle S^{-1}(f), h_k\rangle_{G\times G'})_{k \in  \Lambda}\, 
\|_{\ell_{2,w}}
\\
& = & \|S^{-1}\|^{-1}A\, \|\, 
(\langle f, {{S}^*}^{-1}h_k\rangle_{H\times H'})_{k \in  
\Lambda}\, \|_{\ell_{2,w}} \, ,
\end{eqnarray*}
and (ii) is proved 
with 
$\wt C = \wt B / \wt A =
C\|{S}\|\|{S}^ {-1}\|$.
\end{proof} 


\subsection{Nonlinear  Widths for Frame Pairs} \label{introfw}


The aim of this paper is to study the asymptotic behavior of specific
nonlinear approximation schemes based on frames 
and to compare them with
other well-known widths. 
Especially, we want to prove frame analogues to the results
obtained in \cite{dns1,dns2} for the nonlinear 
widths associated with classes of
Riesz bases. 

Let $({\mathcal F},  {\mathcal G})$ be a frame 
pair for $(H,w)$ in the sense of Definition \ref{defpair} 
and consider specific $n$-term approximations of the form
\begin{equation} \label{bestframe}
\sigma_n \Big(u, ({\mathcal F}, {\mathcal G})\Big)
:= \inf_{|\Lambda| \le n} \, 
\Big\| \, u-\sum_{k \in
    \Lambda} \langle u, h_k\rangle_{H\times H'} \, g_k\, \Big\|_H.
\end{equation}
We do not allow arbitrary expansions in terms of 
the $g_k$ involving at most $n$ nonvanishing coefficients. 
The reason is that, for practical applications, only these 
canonical representations are used. 
Furthermore, to end up with a reasonable notion of a width we need 
to restrict us to stable frame pairs.
\\
In what follows we shall use the conventions: if $F$ is a subspace of $G$
and if $S: \, G \to H$ is an isomorphism then
we equip the subspace $S(F)$ with the quasi-norm 
$\| \, S(f)  \, |S(F)\|:= \|\, f\, |F\|$.
Furthermore, if $K $ is a subspace of $S(F)$ we endow it with the quasi-norm
of $S(F)$.

\begin{definition} \label{defiwidth}
Let $G$ and $H$ be separable Hilbert spaces and 
let $S: ~G \to H$ be an isomorphism.
Let $F$ be a quasi-normed subspace of $G$. For a given constant $C\ge 1$ we denote by 
$\cK_C$ the set of all subspaces $K \subset S(F)$ such that
the inequality
\begin{equation} \label{defkcont}
e_n^{\rm cont} (I, S(F) , H) \leq C
e_n^{\rm cont} (I, K , H)
\end{equation}
holds for all $n$. 
Then, for  $n \in \Nb$, the {\em nonlinear frame width}
$e_{n,C}^{\rm frame}(S, F, H)$ of the operator $S$ is defined by
\begin{equation}
\label{defnlframe}
e_{n,C}^{\rm frame}(S, F, H)
:= \inf \Big\{\sup_{\|f\|_F\leq 1}\, 
\sigma_n \Big(S(f),({\mathcal F}, {\mathcal G})\Big) \mid 
\: (\cF, \cG) \in  {\cal P}_C(K),~ K \in \cK_C\, \Big\}\, .
\end{equation}
\end{definition}

\begin{rem} \label{remarkdrei}
We comment on this definition. To get a reasonable lower bound for 
$e_{n,C}^{\rm frame}(S, F, H)$ we need to restrict ourselves to frame pairs which
are stable with respect to subspaces $K$ of $S(F)$ 
which are not too small.
``Not too small'' is expressed by the inequality {\rm (\ref{defkcont})}.
\end{rem}

In the above definition we decided for 
the manifold widths because they have some 
nice properties.
These  widths $e^\cont_n$   are particular examples of $s$-numbers
in the sense of {\rm Pietsch \cite{Pi74}},  see also {\rm \cite{Ma90}}. 
One of the interesting properties consists in the inequality
\begin{equation}   \label{mult}
e_n^\cont (T_2 \circ T_1 \circ T_0, E_0, F_0) \le \| \, 
T_0\, \|\, \, \|\, T_2\, \|
\, e_n^\cont (T_1, E,F)\, , 
\end{equation} 
where $T_0 \in \cl (E_0,E)$, $T_1 \in \cl (E,F)$, $T_2 \in \cl (F,F_0)$
and $E_0,E,F,F_0$ are arbitrary quasi-Banach spaces. 
As a consequence one obtains that 
the asymptotic behavior of the manifold widths 
remains unchanged under isomorphisms.
A similar result is true in case of our 
nonlinear frame widths.
As a consequence we can concentrate on the 
investigation of identity operators
in what follows.

\begin{lemma}\label{basic} 
Let $G$ and $H$ be separable Hilbert spaces 
and let $S: G \to H$ be an isomorphism.
Let $F$ be a quasi-normed subspace of $G$ and 
let $I:F\to G$ be the identity.  For $C\ge 1$
and
$$
\wt C = C \, ( \Vert S^{-1} \Vert \, \Vert S \Vert )^2 
$$
we obtain 
\begin{equation}\label{e25}
e^{\rm frame}_{n, \wt C} (S,F,H) \le \Vert S \Vert \, 
e^{\rm frame}_{n,C} (I,F,G) 
\end{equation} 
and 
\begin{equation}\label{e25a}
e^{\rm frame}_{n, \wt C} (I,F,G) \le \Vert S^{-1}  \Vert \, 
e^{\rm frame}_{n,C} (S,F,H). 
\end{equation} 
\end{lemma}

\begin{proof} 
We shall prove (\ref{e25}), 
the proof of \eqref{e25a} is very similar. From (\ref{defnlframe}) we can conclude 
that for any   $\varepsilon >0$ we can find a subspace $K \in {\mathcal K}_{C} $
and a frame pair  $(\cF, \cG) \in {\cal P}_C(K)$  for $(G,w)$ such that
\[
\sup_{\|f\|_F \le 1} \inf_{|\Lambda| \le n} \, \Big\| \, f- \sum_{k\in
  \Lambda}\, \langle f, h_k \rangle_{G\times G'} \, g_k\, \Big\|_G 
\le e^{\rm frame}_{n,C}(I,F,G) + \varepsilon \, .
\]
Lemma \ref{framemap} implies  
that  $({S^*}^{-1}(\cF), S({\cG}))$ is a frame pair 
for $(H,w)$ which is contained 
in ${\cal P}_{C_1}(S(K))$, where 
$C_1 = C \, \Vert S^{-1} \Vert \, \Vert S \Vert$. 
We consider the following commutative diagrams:
\[
\begin{CD}
S(F )@> I_1 >> H\\
@VS^{-1}VV @AA S A\\
F @> I_2 >> G
\end{CD}
\qquad \hspace{2cm}\quad
\begin{CD}
K@> I_2 >> G\\
@VSVV @AA S^{-1} A\\
S(K) @> I_1 >> H \, .
\end{CD}
\]
By means of (\ref{mult}) we derive from these diagrams
\[
e^\cont_n (I_1,S(F),H) \le \| S^{-1}\| \,   \| S \| \, 
e^\cont_n (I_2,F,G)
\]
and 
\[  
e^\cont_n (I_2,K,G) \le \| S^{-1}\| \,   \| S \| \, 
e^\cont_n (I_1,S(K),H)\, .
\]
Now our assumption $K \in \cK_C$ yields
\begin{eqnarray*}
e^\cont_n (I_1,S(F),H) &\le &\, \| S^{-1}\| \,   \| S \| \,e^\cont_n (I_2,F,G) \leq  C\, \| S^{-1}\| \,   \| S \| \, 
e^\cont_n (I_2,K,G)\\ 
& \le & C \, (\| S^{-1}\| \,   \| S \|)^2 \, 
e^\cont_n (I_1,S(K),H)\, .
\end{eqnarray*}
In  other words, $S(K)$ belongs to the set $\cK_{ \wt C}$. 
{}From 
\[
\Big\|\, S(f) - \sum_{k \in \Lambda} \langle S(f), {S^*}^{-1}(h_k)
\rangle_{H\times  H'}\, S(g_k)\, \Big\|_H \leq \|S\|\, 
\Big\|\, f-\sum_{k \in \Lambda} \langle f,
h_k\rangle_{g\times G'} \, g_k\, \Big\|_G
\]
it follows that 
\[
e^{\rm frame}_{n, \wt C} (S,F,H) \le \| S\| \,  
e^{\rm frame}_{n,C} (I,F,G) \, . 
\]
\end{proof}

We finish this section by proving two additional properties of nonlinear frame widths that will be used later on in Section
\ref{basic3}. 

\begin{lemma}\label{zusatz}
Let $G_1,G_2, H_1,H_2$ be Hilbert spaces and let 
$S_i \in \cl (F_i,H_i)$, $i=1,2$, be isomorphisms.
Let $F_1,F_2$ be quasi-normed subspaces of $G_1$ and $G_2$, respectively.
Furthermore we suppose  $T_1\in\cl (F_1,F_2) $,  $T_2 \in \cl (H_2,H_1)$ 
and both are isomorphisms.
Finally, we assume that we can decompose  $S_1= T_2 \circ S_2\circ T_1$.
Then,
\begin{equation}\label{eq-60}
e^{\rm frame}_{n, \wt C} (S_1,F_1,H_1) \le  \| T_2\| \, \| T_1 \|  \, 
e^{\rm frame}_{n,C} (S_2,F_2,H_2)
\end{equation} 
 holds with $\wt C = C \, \Vert T_2^{-1} \Vert \, \Vert T_2 \Vert$.
\end{lemma}

\begin{proof}
Corresponding to our assumptions we have
the following commutative diagram:
\[
\begin{CD}
F_1 @>> S_1 > H_1\\
@VT_1VV @AA T_2 A\\
F_2 @>> S_2 > H_2\, .
\end{CD}
\]
By definition, for any   $\varepsilon >0$ we can find a subspace $K \in {\mathcal K}_{C} \subset G$
and a frame pair  $(\cF, \cG) \in {\cal P}_C(K)$  for $(H_2,w)$ such that
\[
\sup_{\|f\|_{F_2} \le 1} \inf_{|\Lambda| \le n} \, \Big\| \, S_2 f- \sum_{k\in
  \Lambda}\, \langle S_2 f, h_k \rangle_{H_2\times H_2'} \, g_k\, \Big\|_{H_2} 
\le e^{\rm frame}_{n,C}(S_2,F_2,H_2) + \varepsilon \, .
\]
Lemma \ref{framemap} implies  
that  $({T_2^*}^{-1}(\cF), T_2({\cG}))$ is a frame pair 
for $(H_1,w)$ which is contained 
in ${\cal P}_{\wt C}(T_2(K))$, where 
$\wt C = C \, \Vert T_2^{-1} \Vert \, \Vert T_2 \Vert$. 
We put
\[
u_k := {T_2^*}^{-1}(f_k)\qquad \mbox{and} \qquad v_k := T_2(g_k)\, .
\]
Consequently
\begin{eqnarray*}
\Big\| \, S_1 g & - & \sum_{k\in
  \Lambda}\, \langle S_1 g, u_k \rangle_{H_1\times H_1'} \, v_k\, \Big\|_{H_1} 
\\
&\le &  \| T_2\|\, 
\Big\| \, S_2 (T_1 g) - \sum_{k\in
  \Lambda}\, \langle S_2 (T_1g), T_2^* u_k \rangle_{H_2\times H_2'} 
\, T_2^{-1}v_k\, \Big\|_{H_2} 
\\
& \le & \| T_2\| \, (e^{\rm frame}_{n,C}(S_2,F_2,H_2) + \varepsilon) \, ,
\end{eqnarray*}
if $\| \, T_1 g\, \|_{F_2} \le 1$.
A homogeneity argument yields
\[
\sup_{\|g\|_{F_1} \le 1} \inf_{|\Lambda| \le n} \,
\Big\| \, S_1 (g)  -  \sum_{k\in
  \Lambda}\, \langle S_1 g, u_k \rangle_{H_1\times H_1'} \, v_k\, \Big\|_{H_1} 
 \le  \| T_2\| \, \| T_1 \| \, e^{\rm frame}_{n,C}(S_2,F_2,H_2) 
\]
which proves our claim.
\end{proof}

\begin{lemma}\label{zusatz2}
Let $U$ be a closed subspace of the Hilbert space $H$ equipped with the same norm as $H$. Let $G$ be a Hilbert space and let
$S:\, G\to H$ be an isomorphism. If  $F$ is a subset of $S^{-1}(U)$, then
\[
e^{\rm frame}_{n, C} (S,F,U) \le  \, e^{\rm frame}_{n,C} (S,F,H) 
\]
follows.
\end{lemma}

\begin{proof}
The Hilbert space $H$ can be written as the orthogonal sum of $U$ and its orthogonal complement $V$. By $P$ we denote the orthogonal projection onto $U$.
Let $(\cF,\cG)$ be a frame pair for $(H,w)$.
Then the elements $f \in U$ can be written in the form
\[
f = \sum_{k=1}^\infty \langle f , h_k\rangle\, P g_k \, .
\]
The norm equivalences (\ref{normpair}) remain unchanged. Hence,
$(\cF, P(\cG))$ is a frame pair for $(U,w)$ with constants $\wt A, \wt B$
and $A \le \wt A \le \wt B \le B$. Concerning the stability it is enough to notice that
only subsets $K$ of $S(F) \subset U$ come into consideration.
\end{proof}

\section{Main Results} \label{jetzt} 

In this section, we want to state and to prove the main results of this
paper. The first theorem  is a general result  
for arbitrary  Hilbert spaces $H$ and $G$
that  clarifies the relationships of the manifold  widths $ e_{n}^\cont
(S,F,H)$
with  the nonlinear frame widths 
$e^{\rm frame}_{n,C} (S,F,H)$.  The second
theorem deals with the more specific situation 
of function spaces on Lipschitz
domains contained in $\R^ d$ and  provides upper and lower bounds for
$e^{\rm frame}_{n,C} (S,B^{-s+t}_q(L_p(\Omega)),H^{s}_0(\Omega))$.

\begin{thm}\label{t1}        
Let $H$ and $G$ be separable Hilbert spaces.
Let $S: G \to H$ be an isomorphism. 
Suppose that the embedding $F\hookrightarrow G$ is compact. Then 
for all $C \ge 1$ and all $n \in \Nb$, we have 
\begin{equation}     \label{e22}     
e_{4n+1}^\cont  (S,F,H) \le 2 \, C^2  \, 
e^{\rm frame}_{n,C} (S,F,H)\, .
\end{equation}
\end{thm}

\begin{thm} \label{basic1}
Let  $\Omega$ be a bounded Lipschitz domain contained in $\R^d$. Let 
$0 < p, q \le  \infty$, $s >0$,  
and $t> d(\frac{1}{p} - \frac{1}{2})_+$.
Let $S: {H}^{-s} (\Omega) \to {H}^{s}_0 (\Omega)$
be an isomorphism.
Then there exists a number $C^*$ such that for any $C\ge C^*$ we have
\[
e_{n,C}^{\rm frame} (S,{B}^{-s+t}_{q} (L_p (\Omega)), 
{H}^{s}_0 (\Omega) ) 
\asymp  n^{-t/d}\, .
\]
\end{thm}

\begin{rem} \label{nochmehrbla}
\begin{itemize}
\item[(i)]
The number $C^*$ depends on $\Omega$. 
It is known that for any Lipschitz domain there 
exists an appropriate frame pair as it is needed here. 
However, optimal 
estimates about the stability seem to be not known.
\item[(ii)]
For exact definitions of the  distribution spaces defined on Lipschitz domains 
we refer to the Appendix and to {\rm  \cite{dns2}}
\item[(iii)] Theorem {\rm  \ref{basic1}} is
a frame analogue to
Theorem {\rm  4}
in {\rm \cite{dns2}}. 
In {\rm \cite{dns2}}, it has been shown that if the
domain $\Omega$ is chosen in such a way 
that the spaces ${B}^{-s+t}_{q} (L_p
(\Omega))$ and ${H}^{-s} (\Omega)$ allow a discretization by one common
wavelet system $\tilde{\mathcal R}^*$, then also
$$
e_{n,C}^\non (S,{B}^{-s+t}_{q} (L_p (\Omega)), {H}^{s}_0 (\Omega) ) 
\asymp  n^{-t/d}
$$
holds for $C$ sufficiently large.
We see that the  restrictive condition 
on the domain that was needed  in the Riesz 
basis case can be dropped in the frame setting.
\item[(iv)]
Our proof of the upper bounds in Theorem {\rm \ref{basic1}} is constructive.
One may always use the frame pair constructed in 
Lemma {\rm \ref{framepairstandard}}
below.
\end{itemize}
\end{rem}


\section{Examples}
\label{example}

In this section, we apply the analysis presented above to two classical examples, i.e., the Poisson equation in a Lipschitz domain and 
the single layer potential equation on the unit circle.

\subsection{The Poisson Equation}


We consider the Poisson  equation in a bounded 
Lipschitz domain $\Omega$ contained in $\R^d$
\begin{eqnarray}
-\triangle u &=&f \quad  \mbox{in}\quad \Omega \label{Poisson}\\
           u&=&0 \quad \mbox{on} \quad \partial \Omega. \nonumber
\end{eqnarray}

\noindent
As usual, we study  (\ref{Poisson}) in the weak formulation. Then, it can be
shown that  the operator ${\cal A}= \triangle:~H_0^1\longrightarrow H^{-1}$ is
boundedly invertible, see, e.g., \cite{hack}  for details.
Hence Theorem \ref{basic1} applies with $s=1$, so that 
\[
e_{n,C}^{{\rm frame}} (S, B^{-1+t}_q(L_p (\Omega)), H^1_0(\Omega)) 
\asymp n^{-t/d}\, 
\]
if $t> d\, (\frac 1p - \frac 12)_+$.


\subsection{The Single Layer Potential}


As a second example we shall deal with an integral equation.
Let $\Gamma$ be the unit circle.
Then we consider the Fredholm integral equation of the first kind
\[
\A f (x):=  - \frac{1}{2\pi}\, \int_\Gamma \log |x-y|\, f(y)\, 
d\Gamma_y = \varphi (x)\, , \qquad x \in \Gamma\, .
\]
The left-hand side is called the {\em single layer potential.}
The following is known, cf. e.g. \cite{Co88}: the operator 
$\A$ belongs to $\cl (H^{-1/2}(\Gamma),H^{1/2}(\Gamma))$,
where $H^{1/2}(\Gamma)$ is the collection of all functions $g \in L_2 (\Gamma)$
such that
\[
\int_{\Gamma} \int_{\Gamma} \frac{|g(x)-g(y)|^2}{|x-y|^2}\, d\Gamma_x \, d\Gamma_y
<\infty
\]
and $H^{-1/2}(\Gamma)$ its dual.
Furthermore, $\A$ is a bijection of $H$ onto $G$ where
\[
G:= \{ g \in H^{1/2}(\Gamma): \: \int_\Gamma g(y)\, d\Gamma_y =0\}
\quad \mbox{and} \quad
H := \{ g \in H^{-1/2}(\Gamma): \: \langle g, 1\rangle =0 \}\,.
\]
The space $G$ can be interpreted as the quotient space 
$H^{1/2}(\Gamma)/\R$ of 
$H^{1/2}(\Gamma)$ with $\R$ (the constants)
and $H$ can be interpreted as the quotient space 
$H^{-1/2}(\Gamma)/\R$ of 
$H^{-1/2}(\Gamma)$ with $\R$.
By $S$ we denote $\A^{-1}$, defined on $G$ with values in $H$.
Now we investigate 
$e_{n,C}^{\rm frame} (S,F,G)$ where $F$
is chosen to be the quotient space of the Besov space
$ B^{t+ 1/2}_{q}(L_p (\Gamma))$ and the constants, see Subsection
\ref{circle} for a definition of $ B^{t+ 1/2}_{q}(L_p (\Gamma))$.
We put
\[
 Y^{s}_{q}(L_p (\Gamma)) := \{ g \in  B^{s}_{q}(L_p (\Gamma)): 
\: \langle g, 1\rangle_\Gamma =0 \}\,.
\]
The same principles as above apply.
Again we use a commutative diagram
\begin{eqnarray}\label{final}
H^{1/2}(\Gamma)/\R & \stackrel{\hbox{$S$}}{\longrightarrow} & H^{-1/2}(\Gamma)/\R 
\nonumber
\\
I \,  & \nwarrow \qquad  \nearrow & S_F \\
& F:= Y^{t+ 1/2}_{q}(L_p (\Gamma)). &
\nonumber
\end{eqnarray}
Here $I$ denotes the identity and $S_F$ the restriction of 
$S$ to $F$. Then the outcome is as follows.

\begin{thm}
\label{basic2}
Let $0 < p, q \le  \infty$ and $t> (\frac{1}{p} - \frac{1}{2})_+$.
Then there exists a number $C^*$ such that for any $C\ge C^*$ we have
\[
e_{n,C}^{\rm frame} (S,{Y}^{t+1/2}_{q} (L_p (\Gamma)), H) 
\asymp  n^{-t}\, .
\]
\end{thm}

\begin{rem}
There are far-reaching extensions concerning the theory 
of the mapping properties of the single layer potentials. In particular, 
much more general curves and surfaces are discussed.
We refer to {\rm \cite[Sect.~20]{T01}} for the discussion of these properties 
in the framework of $d$-sets. 
\end{rem}


\section{Proofs}\label{Beweis}


\subsection{Proof of Theorem \ref{t1}} \label{pt1}

First we deal with Theorem 1. 
Here we shall work in the framework of Hilbert frame pairs. 
Hence we consider sequences $(g_k)_k$ and $(h_k)_k$ 
in a (separable) Hilbert space $H$ such that 
\begin{equation}  \label{lo01}
f = \sum_{k=1}^\infty (f, h_k) g_k 
\end{equation} 
for all $f \in H$, compare with Remark \ref{ganzviellaber} (ii).  
By (17) 
we may assume that 
\begin{equation}  \label{lo02}
\left\Vert \sum_{k=1}^\infty c_k g_k \right\Vert^2 \le B^2 \cdot 
\sum_{k=1}^\infty  c_k^2 
\end{equation} 
for arbitrary $(c_k)_ {k \in \Nb} \in \ell_2(\Nb)$. 
Moreover, we assume that  the representation 
\eqref{lo01} is stable on $K \subset H$ in the sense that 
\begin{equation}  \label{lo03}
A^2 \sum_{k \in \Lambda} |(f,h_k)|^2 \le 
\left\Vert \sum_{ k\in \Lambda}  (f, h_k) g_k \right\Vert^2
\end{equation} 
for arbitrary $f \in K$ and $\Lambda \subset \Nb$.
Moreover we assume that 
\begin{equation}  \label{lo04}
\frac{B}{A} \le C .
\end{equation} 
We consider particular $n$-term approximations  
of $f \in K$ by subsums of \eqref{lo01} and their 
error
\begin{equation}  \label{lo05}
\sigma_n (f) = \inf_{|\Lambda| \le n}\left\Vert f - \sum_{k \in \Lambda} 
(f, h_k) g_k \right\Vert .
\end{equation} 
We define 
\begin{equation}  \label{lo06}
e_{n,C} (K,H) = \inf_{(g_k)_k, (h_k)_k} 
\sup_{f \in K}  \sigma_n (f) ,
\end{equation} 
with the understanding that \eqref{lo01}-\eqref{lo05} 
hold true.  Moreover, we define
\begin{equation} \label{dreiundvierzig}
e_{n}^{\rm cont} (K, H):= \inf_{N_n, \varphi_n} \sup_{u  \in K} \|\varphi_n(N_n(u))-u\|,
\end{equation}
where the infimum runs over all continuous mappings $\varphi_n~:~\R^n\rightarrow H$ and $N_n~:~K \rightarrow \R^n.$ 
Then the following result is a frame analogue of 
Proposition 1 from \cite{dns2}. 

\begin{proposition} \label{prop1} 
Assume that $K \subset H$ is compact and $C \ge 1$. 
Then 
\begin{equation}  \label{lo07}
e_{4n+1}^{\rm cont} (K, H) \le 2 C e_{n,C} (K,H) .
\end{equation} 
\end{proposition} 

\begin{proof} 
Assume that $K$, $n$, and $C \ge 1$ are given.
Let $\e >0$. 
Then there exist sequences $(g_k)_k$ and $(h_k)_k$ 
in $H$ such that \eqref{lo01}-\eqref{lo04} as well as 
\begin{equation}  \label{lo08}
\sup_{f \in K} 
\inf_{|\Lambda| \le n} \Vert f - \sum_{ k \in \Lambda} 
(f, h_k) g_k \Vert \le e_{n, C} (K, H) +\e 
\end{equation} 
hold. 
Since we only consider $f \in K$, we can always assume 
that the index set $\Lambda$ is a subset of 
$\{ 1, 2, \dots , N \}$. We only loose another $\e$. 
Here $N$ might be large, but is finite. 
We write 
\begin{equation}  \label{lo09}
L_N(f) = \sum_{k=1}^N (f, h_k) g_k 
\end{equation} 
and obtain 
\begin{equation}  \label{lo10}
\sup_{f \in K} \Vert f - L_N(f) \Vert \le \e 
\end{equation} 
and
\begin{equation}  \label{lo11}
\sup_{f \in K} \inf_{|\Lambda | \le n}  \left\Vert L_N(f) - 
\sum_{k \in \Lambda} (f, h_k) g_k \right\Vert  \le 
e_{n,C} (K, H) + 4 \e .
\end{equation} 
For the $n$-term approximation in \eqref{lo11} we also 
write
\begin{equation}  \label{lo12}
f_n^* = \sum_{k \in \Lambda} a_k g_k, 
\end{equation} 
hence $a_k = (f,h_k)$ and $|\Lambda | =n$ for each $f \in K$ 
and 
\begin{equation}  \label{lo13}
\sup_{f \in K}  \Vert L_N(f) - 
f_n^*   \Vert  \le 
e_{n,C} (K, H) + 4 \e .
\end{equation} 
For the proof we may assume that $A=1$. 
We consider the modification $L_N^*$ of $L_N$ defined by
\begin{equation}  \label{lo14}
L_N^* (f) = \sum_{k=1}^N a_k^* g_k , 
\end{equation} 
where 
$a_k^* =a_k$ if $|a_k | \ge 2  \b$ and 
$a_k^* =0$ if $|a_k| \le \b$. 
To obtain a continuous dependence of $a_k^* $ from 
$a_k$ and, hence, a continuous mapping 
$L_N^* : H \to H$,  we define
$$
a_k^* = 2 \sgn a_k \cdot ( |a_k| - \b) 
$$
if $|a_k| \in (\b, 2\b)$. 
The number $\b >0$ will be  defined later. 

Assume that for $f \in K$ there are $m>n$ of the $a_k$
with $|a_k| \ge \beta$. 
Then 
$$
L_Nf - f_n^* = \sum_{k \in \tilde \Lambda} a_k g_k ,
$$
where $\tilde \Lambda$ contains at least $m-n$ elements 
with  $|a_k| \ge \b$. 
Then we obtain from \eqref{lo03} 
$$
\Vert L_Nf -f_n^* \Vert \ge (m-n)^{1/2} \beta 
$$
and with \eqref{lo13} we get 
\begin{equation}  \label{lo15} 
m-n \le \frac{1}{\beta^2} ( e_{n,C}(K, H) +4\e)^2 .
\end{equation} 
Now we consider the sum $\sum_{|a_k|<\beta} a_k^2$ for $f \in K$. 
We distinguish between those $k$ that are used for $f_n^*$ 
(there are at most $n$ of those $k$) and the other indices and 
obtain 
\begin{equation}   \label{lo16} 
\sum_{|a_k|<\beta} a_k^2 \le n \beta^2 + (e_{n,C}(K,H) +4\e)^2 .
\end{equation} 
Now we are ready to estimate $\Vert L_N^*(f) - L_N(f)\Vert$ 
for $f \in K$. Observe that $|a_k^* - a_k| \le \beta$ 
for any $k$. We obtain 
$$
\Vert L^*_N (f) - L_N(f) \Vert \le 
B ( m\beta^2 + n \beta^2 + (e_{n,C}(K,H) +4\e)^2 )^{1/2} .
$$
Using the estimate \eqref{lo15} for $m$,  we obtain 
$$
\Vert L^*_N (f) - L_N(f) \Vert \le 
B ( 2  n \beta^2 + 2 (e_{n,C}(K,H) + 4 \e)^2)^{1/2} . 
$$
Now we define $\beta$ by
$$
n\beta^2 = ( e_{n,C}(K,H) +4\e)^2 
$$
and obtain the final error estimate (where we replace, for general $A$,
the number $B$ by $B/A$) 
$$
\Vert L_N^* (f) - L_N(f) \Vert \le 
\frac{2B}{A} \,  (e_{n,C}(K,H) +4 \e) .
$$
In addition we obtain 
$$
m \le 2n 
$$
and therefore $L_N^*$ yields a
continuous $2n$-term approximation of 
$f \in K$ with error at most
$$
\sup_{f \in K} \Vert L_N^* (f) - f \Vert \le 
\frac{2B}{A} \,  (e_{n,C}(K,H) + 4\e) + \e .
$$
The mapping $L_N^*$ is continuous and the image is 
a complex of dimension $2n$, see, e.g.,  \cite{DKLT93}. 
Hence we have an upper bound for the so-called Aleksandrov 
widths, see \cite{DKLT93} and \cite{St74}.  
By the famous theorem of N\"obeling, 
any such mapping can be factorized as
$L_N^*=\phi_{4n+1} \circ N_{4n+1}$ where 
$N_{4n+1} : K \to \R^{4n+1} $ and 
$\phi_{4n+1} : \R^{4n+1} \to H$ are continuous. 
Hence the result is proved. 
\end{proof} 

\subsubsection*{Proof of Theorem \ref{t1}} \label{pt1a}

First we observe that 
$$
e_{4n+1}^\cont (S,F, H) = e_{4n+1}^\cont (I, S(F), H) .
$$
Condition  \eqref{defkcont}  implies that 
$$
e_n^\cont (I, S(F), H) \le C \, e_n^\cont 
(I, K, H)  = C \, e^\cont_n (K,H) 
$$
so that Proposition~\ref{prop1} yields
$$
e_{4n+1}^\cont (K,H) \le 2C \, e_{n,C}(K,H) 
\le 2C \, e_{n,C}^\frame (I, S(F), H) .
$$
We also have
$$
e_{n,C}^\frame (I, S(F), H) = 
e_{n,C}^\frame (S, F, H), 
$$
hence we finally  obtain 
$$
e_{4n+1}^\cont (S,F, H) \le 2C^2 \,
e_{n,C}^\frame (S, F, H). 
$$

\subsection{Proof of Theorem \ref{basic1}} \label{pbasic1}

We want to make a  general remark concerning  the notation in advance. 
In what follows we will use the symbol
$\langle \cdot, \cdot \rangle$ for different duality pairing. 
Which one will be always clear from the context. So we avoid indices.

\subsubsection{Lower Bounds} \label{lbasic}

The proof of the lower bound follows by 
combining Theorem \ref{t1}  with the  following proposition proved in \cite{dns2}, see also \cite{DHM89,DKLT93,DD96}:

\begin{proposition}\label{t4} 
Let $\Omega \subset \Rd$ be a bounded Lipschitz domain.
Let $0 < p, q \le \infty$, $s>0$, and 
$$
t> d \left(  \frac{1}{p} - \frac{1}{2} \right)_+ \, .
$$
Then 
\[
e_n^\cont (S,{B}^{-s+t}_{q} (L_p (\Omega)), {H}^{s}_0 (\Omega) ) \asymp
n^{-t/d}\, .
\]
\end{proposition}

\subsubsection{Upper Bounds} \label{ubasic}

The proof of the upper bound  turns out to be a 
little bit more complicated.
However, let us mention that our proof is constructive.
As a first step we reduce the proof 
of Theorem \ref{basic1} to the proof of 
the following

\begin{thm}   \label{th5} 
Let $\Omega$ be as above.
Let $0 < p, q \le \infty$, $s \in \R $ 
and suppose that  
\[
t> d\, \Big(\frac{1}{p} - \frac 12 \Big)_+
\] 
holds. Then there exists a number $C^*$ 
such that for any  $C\ge C^*$  we have 
\[
e_{n,C}^{\rm frame} (I, B^{s+t}_{q} 
(L_p(\Omega)), B^s_{2} (L_2(\Omega)))   \lsim
n^{-t/d}\,  . 
\]
\end{thm} 

\noindent
{\bf Proof of Theorem \ref{basic1}}. Since  $H^{-s} (\Omega) = B^{-s}_2 (L_2(\Omega))$, 
cf. Remark \ref{tilde}, 
Theorem \ref{th5} yields that 
\[
e_{n,C}^{\rm frame} (I,  B^{-s+t}_{q} (L_p (\Omega)),   
H^{-s}(\Omega) )    
\lsim  n^{-t/d}. 
\]
Since $S: H^{-s} (\Omega) \to H^s_0 (\Omega) $ is an isomorphism,  
Lemma \ref{basic} implies the desired result. 
{\eproof}


\subsubsection{Widths and Discrete Besov Spaces}
\label{discr}


The proof of Theorem \ref{th5} requires several preparations.
First of all, let  us fix some notation. 
Let $0< p,q \le \infty$ and let $s \in \R$. 
Let $\nabla := (\nabla_j)_{j=-1}^\infty$ be a sequence 
of subsets of finite cardinality of
the set $\{1,2,\ldots \,  , 2^d-1\} \times \Z^d$. 
We suppose that there exist $0< C_1 \le C_2$ 
and $J \in \Nb$ such that the 
cardinality $|\nabla_j|$ of $\nabla_j$ satisfies
\begin{equation}\label{301}
C_1 \le  2^{-jd} \, |\nabla_j| \le C_2 
\qquad \mbox{for all}\quad j \ge J \, . 
\end{equation} 
Then $b^s_{p,q}(\nabla)$,  where $ 0< q < \infty$, 
denotes the collection of all sequences
$a = (a_{j,\lambda})_{j,\lambda}$ of complex numbers such that
\begin{equation}\label{300}
\| \, a \, \|_{b^s_{p,q}}  := \left(  \sum_{j=-1}^\infty  
2^{{j(s+d( 1/2 -1/p))q}} 
\bigg( \sum_{\lambda \in \nabla_j}
| \, a_{j,\lambda} |^p \bigg)^{q/p}
\right)^{1/q} <\infty\,.
\end{equation}
For $q=\infty$, we  use the usual modification
\begin{equation} \label{300a}
\| \, a \, \|_{b^s_{p,\infty}}:= \sup_{j=-1,0,1,\ldots}
2^{{j(s+d({1}/{2}-{1}/{p}))}} 
\left(\sum_{\lambda \in \nabla_j}|a_{j,\lambda}|^p\right)^{1/p}
<\infty.
\end{equation}
In our paper \cite{dns2} we have dealt with 
several types of  widths of embeddings of those 
discrete Besov spaces. A few of the results we obtained there will 
be recalled now.  

\begin{proposition}   \label{discrete} 
Let $0 < p, q \le \infty$
and $s \in \R $.  Suppose that
\begin{equation}\label{00}
t> d\, \Big(\frac{1}{p} - \frac 12 \Big)_+\, .
\end{equation}
It holds
\[
e_{n}^\cont (I, b^{s+t}_{p,q}({\nabla}),  
 b^s_{2,2}({\nabla})) \asymp
e_{n}^\non (I, b^{s+t}_{p,q}({\nabla}),  
 b^s_{2,2}({\nabla}))   \asymp n^{-t/d} \, .
\]
\end{proposition} 

\begin{rem}
Of course, the constants in the above inequalities depend on $\nabla$
(and therefore on $C_1,C_2$ and $J$) as well as on $s,t,p$ and $q$.
But this will play no role in what follows.
\end{rem}


\subsubsection{Frame Pairs for Sobolev Spaces on Domains}


Now we turn to the construction of frame pairs 
for Sobolev spaces with some additional features.

Let $s \in \R$ be fixed and 
let 
\begin{equation}
\label{psi}
\Psi := \Big\{\varphi_{k}, \widetilde{\varphi}_{k}: \: 
k \in \Zd\Big\}\, 
\cup\,  \Big\{ \psi_{i,j,k}, \widetilde{\psi}_{i,j,k}: 
\quad i =1, \ldots 2^d-1, \, j=0,1,2\ldots\, , k \in \Zd\Big\}\, , 
\end{equation}
be a biorthogonal wavelet system such that 
the parameter $r$, controlling the smoothness and the 
moment conditions, satisfies $r>|s|$, see Proposition
\ref{wavelets} in the Appendix. 
Here, as always in this subsection we shall 
use $H^s(\Omega)= B^s_2(L_2(\Omega))$
in the sense of equivalent norms, see the Appendix. 
We suppose
\[
\supp \varphi\, , \, \supp \psi_i\, , \,  
\supp \widetilde\varphi\, , \, 
\supp \widetilde{\psi}_i
\subset [-N,N]^d\, , \qquad i=1, \ldots \, 
2^d-1\, .
\]
By $B(x^0,R)$ we denote a ball with radius $R$ and center $x^0$.
We may assume
$\Omega \subset B(x^0,R)$ for some $R>0$ and $x^0 \in \Omega$. 
Rychkov {\rm \cite{Ry}} has proved that in case of a bounded Lipschitz domain there exists a linear and continuous extension operator $\ce \in \cl (H^s (\Omega) \to H^s (\Rd))$. In addition we may assume that 
\begin{equation}\label{eq105}
\supp \, \, \ce f \, \subset B(x^0,2R)
\end{equation}
holds for all $f\in H^s(\Omega)$.
Now we turn to the wavelet decomposition of $\ce f$. 
Defining
\[
\Lambda_j  :=  \Big\{ k \in \Zd: \quad |2^{-j}k_i - x_i^0|\le  2R+ 2^{-j}N
\, , \quad i=1, \ldots \, , d\Big\}\, , \qquad j=0,1, \ldots \, ,
\
\]
we obtain for given $f \in H^s (\Omega)$ 
\begin{equation}\label{eq201}
\ce f = \sum_{k \in \Lambda_0} \langle \ce f, 
\widetilde{\varphi}_k \rangle \, \varphi_k
 + \sum_{i=1}^{2^d-1} \, \sum_{j=0}^\infty 
\sum_{k \in \Lambda_j} \langle \ce f, 
\widetilde{\psi}_{i,j,k} \rangle\, \psi_{i,j,k} 
\qquad (\mbox{convergence in} \quad \S')
\end{equation}
and 
\begin{eqnarray}\label{eq202}
\| \, \ce f\, |H^s (\Rd)\| & \asymp &  \Big(\sum_{k \in \Lambda_0} 
| \langle \ce f, \widetilde{\varphi}_k \rangle |^2\Big)^{1/2}  +
\\
&&
\bigg(\sum_{i=1}^{2^d-1} \, \sum_{j=0}^\infty \, 
2^{2js} 
\Big(\sum_{k \in \Lambda_j} | \langle \ce f, 
\widetilde{\psi}_{i,j,k} \rangle |^2\Big)\bigg)^{1/2} < \infty\, .
\nonumber
\end{eqnarray}
This can be rewritten by using
\begin{eqnarray}\label{nabla1}
\nabla_{-1} & := & \Lambda_0 \\
\label{nabla2}
\nabla_j & := & \Big\{
(i,k) : \quad 1 \le i \le 2^d-1\, , \quad k \in \Lambda_j
\Big\}\, , \quad j=0,1,\ldots \, , 
\end{eqnarray}
$\psi_{j,\lambda}  := \psi_{i,j,k}$, if $\lambda = (i,k) \in \nabla_j$,
$j \in \Nb_0$, and 
$\psi_{j,\lambda} 
:= \varphi_k $ if $\lambda =k \in \nabla_{-1}$. Similarly in
case of the dual basis.
Then (\ref{eq201}), (\ref{eq202}) read as
\begin{equation}\label{eq203d}
\ce f = \sum_{j=-1}^\infty 
\sum_{\lambda \in \nabla_j} \langle \ce f, 
\widetilde{\psi}_{j,\lambda} \rangle\, \psi_{j,\lambda} 
\qquad (\mbox{convergence in} \quad \S')
\end{equation}
and 
\begin{equation}\label{eq204}
\| \,  f \, |H^s(\Omega)\| \asymp \| \, \ce f \, |H^s(\Rd)\|  \asymp \| \, ( \langle \ce f, \widetilde{\psi}_{j,\lambda} \rangle )_{j,\lambda}\, \|_{b^{s}_{2,2} (\nabla)} \, .
\end{equation}
Let $\cx_\Omega$ denote the characteristic function of $\Omega$.
We put
\begin{equation}\label{gframe}
g_{j,\lambda}:= \cx_\Omega \, \psi_{j,\lambda}\, ,  \qquad  j=-1,0,1, \ldots \,, 
\quad  \lambda \in \nabla_j\, .
\end{equation}
For  $M \in \Nb$ we have
\[
\sum_{j=-1}^M \sum_{\lambda \in \nabla_j} \langle \ce f, 
\widetilde{\psi}_{j,\lambda} \rangle\,  g_{j,\lambda}   = 
\Big(\sum_{j=-1}^M \sum_{\lambda \in \nabla_j} \langle \ce f, 
\widetilde{\psi}_{j,\lambda} \rangle\, \psi_{j,\lambda} \Big)_{{\Big|}_{\Omega}}
\]
and consequently 
\[
\lim_{M\to \infty} \sum_{j=-1}^M \sum_{\lambda \in \nabla_j} \langle \ce f, 
\widetilde{\psi}_{j,\lambda} \rangle\,  g_{j,\lambda}   = (\ce f)_{|_\Omega} = f
\]
in $H^s (\Omega)$.
Let $\ce^*$ denote the adjoint of $\ce$. 
Define 
\begin{equation}\label{hframe}
h_{j,\lambda} = \ce^* (\widetilde{\psi}_{j,\lambda}) \, ,  \qquad  j=-1,0,1, \ldots \,, 
\quad  \lambda \in \nabla_j\, .
\end{equation}
Then, taking into account the norm equivalences (\ref{eq204}), it follows that 
$(\cF, \cG)$ satisfies (\ref{atom}) and (\ref{normpair}) for $(H^s (\Omega), b^s_{2,2}(\nabla))$, where
\begin{eqnarray}
\cF & = & \{ h_{j,\lambda} : \quad j=-1,0,1, \ldots \, , \quad 
\lambda \in \nabla_j \} \qquad \mbox{and}\qquad \\
\cG & = & \{g_{j,\lambda} : \quad j=-1,0,1, \ldots \, , \quad 
\lambda \in \nabla_j \}\, .
\label{framepsi}
\end{eqnarray}
Instead of writing $(H,w)$  we used here the notation $(H,\ell_{2,w})$, see
Definition \ref{defpair}.  To obtain a frame pair,  it remains  to establish a suitable reconstruction operator. Due to the norm equivalences  stated in  (\ref{eq202}) and Proposition \ref{wavelets},
it is clear that such an operator $R~:~\ell_{2,w}\longrightarrow  H^s({\mathbb R}^d)$ exists on all of ${\mathbb
R}^d$. Therefore 
$$\tilde{R}~:~b^s_{2,2}(\nabla)\longrightarrow H^s(\Omega), \qquad  a=(a_{j,\lambda})_{(j,\lambda) \in \nabla}
\longmapsto \chi_{\Omega}R(a)$$
does the job. 
We collect our findings in the following lemma.

\begin{lemma} \label{framepairstandard} 
Let $\Omega \subset \R^d$  be a
bounded Lipschitz domain.
Let $\Psi$ be a wavelet system, see {\rm  (\ref{psi})}, such that $r > |s|$, see
Proposition {\rm \ref{wavelets}}.
Let $\cF$ and $\cG$ be defined as in {\rm (\ref{gframe})-(\ref{framepsi})}.
Then $(\cF,\cG)$ is a frame pair  for   $(H^s (\Omega), b^s_{2,2}(\nabla))$,
where $\nabla = \nabla (\Omega)$ is defined in {\rm  (\ref{nabla1}), (\ref{nabla2})}.
\end{lemma}


\subsubsection{Stability of Frame Pairs}


Next we need to investigate the stability of this frame pair
constructed in the previous subsection.
The symbol $\nabla$ will always refer to 
$\nabla = \nabla (\Omega)$ defined in (\ref{nabla1}), (\ref{nabla2}).
Let $0 <p,q\le \infty $ and suppose $t> d (\frac 1p - \frac 12)_+$.
Furthermore, we require that the parameter $r$ of the wavelet system
satisfies
\begin{equation}\label{eq-300}
r > \max\Big(s+t, d\max (0, \frac1p - 1) - s, d \max (0, \frac1p - 1) - 
(s+t)\Big)\, ,
\end{equation}
see Proposition \ref{wavelets}.
We choose a rectangular subset
$\Box$  of   $\Omega$  such that $\dist(\Box, \partial \Omega) >0$. 
Then we define
\begin{equation}\label{nabla3}
\nabla_j^*  :=  \Big\{(i,k) \in \Lambda_j : \quad \supp \psi_{j,\lambda}
\subset \Box \Big\}\, , \quad j=0,1,\ldots \, , 
\end{equation}
Of course, it may happen that $\nabla_j^* = \emptyset$ if $j$ is small.
Let $J \in \Nb$ be a number such that $\nabla_j^* \neq \emptyset$ for all $j \ge J$.
Then we put
\begin{equation}\label{Box}
K  :=  \Big\{ f \in {\cD}' (\Omega): \quad \mbox{there exists}~  (a_{j,\lambda})_{j,\lambda}
\in b^{s+t}_{p,q}(\nabla^*)\: \mbox{s.t.}\quad f = \sum_{j=J}^\infty
\sum_{\lambda \in \nabla_j^*} a_{j,\lambda}\, \psi_{j,\lambda}
\Big\}\, .
\end{equation}
Because of $\dist(\Box, \partial \Omega) >0$ we can extend $f$ by zero outside of $\Omega$ and obtain from Proposition \ref{wavelets}
that $K \subset B^{s+t}_q (L_p (\Omega))$. Again making use of 
Proposition \ref{wavelets} we find that   
\[
\Big\| \, \sum_{(j,\lambda) \in \Lambda} a_{j,\lambda}\, \psi_{j,\lambda}\, 
\Big\|_{H^s(\Omega)} \asymp 
\Big\| \, \sum_{(j,\lambda) \in \Lambda} a_{j,\lambda}\, \psi_{j,\lambda}\, 
\Big\|_{H^s(\Rd)} \asymp 
\| \, (a_{j,\lambda})_{(j,\lambda) \in \Lambda}\, \|_{b^s_{2,2}(\nabla^*)}\, ,
\] 
if $\Lambda \subset \bigcup_{j=J}^\infty \nabla_j^*$. 
Here the constants do not depend on $\Lambda$.

Finally we have to show that $K$ is sufficiently large or more exactly, that
$K\in {\mathcal K}_C$ for some sufficiently large $C$. 
By definition of $K$ the mapping
\[
T : f \mapsto  (\langle f, \widetilde{\psi}_{j,\lambda} 
\rangle)_{(j,\lambda) \in \nabla_j^*}
\]
belongs to $\cl (K, b^{s+t}_{p,q}(\nabla^*))$. Moreover, it is invertible and 
$T^{-1} \in \cl (b^{s+t}_{p,q}(\nabla^*),K)$. 
Once again we shall use the extension operator $\ce$. In addition we apply 
the fact that $\ce$ may be chosen such that
$\ce \in \cl (B^{s+t}_q(L_p (\Omega)),B^{s+t}_q(L_p (\Rd)))$,
cf. Ryshkov \cite{Ry}.
Now we  extend $T$ by defining
\[
T : f \mapsto  (\langle \ce f, \widetilde{\psi}_{j,\lambda} 
\rangle)_{(j,\lambda)  \in \nabla_j}\, .
\]
This extension is again bounded, cf. Proposition \ref{wavelets}.
Let us have a look at the commutative diagram
\[
\begin{CD}
b^{s+t}_{p,q} (\nabla^*) @ > I_1 >> b^s_{2,2}(\nabla)\\
@V T^{-1}VV @AA T A\\
K @> I_2 >> B^s_2 (L_2(\Omega)) \, .
\end{CD}
\]
Because of $\nabla_j^* \subset \nabla_j$, $j \ge J$, there is a natural embedding operator between these sequence spaces, here denoted by $I_1$.
Since $T \in \cl (B^s_2(L_2(\Omega)),b^s_{2,2}(\nabla))$ we can apply (\ref{mult})
and conclude
\begin{equation}   \label{mult1}
e_n^\cont (I_1, b^{s+t}_{p,q} (\nabla^*), b^s_{2,2}(\nabla)) \le \| \, T^{-1}\, \|\, \, \|\, T\, \| \, e_n^\cont  (I_2, K, B^s_2(L_2(\Omega)))\, . 
\end{equation} 
Furthermore
\[
e_n^\cont (I_1, b^{s+t}_{p,q} (\nabla^*), b^s_{2,2}(\nabla^*)) =
e_n^\cont (I_1, b^{s+t}_{p,q} (\nabla^*), b^s_{2,2}(\nabla)) \, .
\]
To explain this we split  $b^s_{2,2}(\nabla))$ 
into $b^s_{2,2}(\nabla^*)$ and its orthogonal complement $U$.
Then the claimed identity  follows from the observation that 
optimal approximations $S_n =\varphi_n \circ N_n$, see (\ref{manifold}), of elements
of $b^{s+t}_{p,q}(\nabla^*))$ are obtained with 
$\varphi_n : \, \R^n \to b^s_{2,2}(\nabla^*)$.
The behavior of the left-hand side in (\ref{mult1}) is known, see Proposition
\ref{discrete}.
As a consequence we obtain
\begin{eqnarray}   \label{mult2}
c_1\, n^{-t/d} & \le &
e_n^\cont (I_1, b^{s+t}_{p,q} (\nabla^*), b^s_{2,2}(\nabla^*)) =
e_n^\cont (I_1, b^{s+t}_{p,q} (\nabla^*), b^s_{2,2}(\nabla)) 
\nonumber
\\
& \le & 
 c_2 \, e_n^\cont  (I_2, K, B^s_2(L_2(\Omega)))
\end{eqnarray} 
with some positive $c_1,c_2$.
Summarizing we have proved that the frame pair $(\cF, \cG)$ from Lemma
\ref{framepairstandard} is admissible in the sense of Definition \ref{defiwidth}
for $C$ sufficiently large.

\begin{lemma} \label{framepairstabil} 
Let $\Omega \subset \R^d$  be a bounded Lipschitz domain. 
Let $\Box $ be a rectangular subset of $\Omega$ such that 
$\dist (\Box, \partial \Omega)>0$. 
Let $s \in \R$, $0< p,q\le \infty$ and
$t> d(\frac 1p - \frac 12)_+$.
Let $\Psi$ be a wavelet system, see {\rm (\ref{psi})}, such that $r$
satisfies {\rm (\ref{eq-300})}, see Proposition {\rm \ref{wavelets}}.
Let $\cF$ and $\cG$ be defined as in {\rm (\ref{gframe})-(\ref{framepsi})}.
Then the frame pair  $(\cF,\cG)$  
is stable with respect 
to the set $K$ defined in {\rm (\ref{Box})}, i.e.
it belongs to ${\cal P}_C (K)$, and it also belongs to $\cK_C \subset 
B^{s+t}_q (L_p (\Omega))$ if $C$ is sufficiently large.
\end{lemma}

\subsubsection{Proof of Theorem \ref{th5}}

To prove Theorem \ref{th5} we shall use the frame pair from  Lemmata
 \ref{framepairstandard} and \ref{framepairstabil}.\\
Let $\Lambda \subset{\nabla}$ be a set of cardinality $n$.
Then 
\begin{eqnarray*} 
\sigma_n (f,(\cF, \cG))_{B^s_2(L_2(\Omega))} &\le & \Big\|
\sum_{(j,\lambda)\not\in  \Lambda} \langle f, {\cal  E}^*\tilde{\psi}_{j,\lambda}
  \rangle\, g_{j,\lambda}\, \Big\|_{B^s_2(L_2(\Omega))} 
\nonumber \\
&\le & c_1 \, \| \, (\langle f, {\cal  E}^*\tilde{\psi}_{j,\lambda}
  \rangle)_{(j,\lambda) \not\in \Lambda}\, \|_{b^s_{2,2}}\, ,
\nonumber
\end{eqnarray*}
where we have once again used (\ref{upperwichtig}). 
By $\O$ we denote the canonical orthonormal basis of $b^0_{2,2}(\nabla)$
and by $e_{j,\lambda}$ its elements, respectively.
For $a \in b_{2,2}^s (\nabla)$ we put
\[
\sigma_n \big( a, \O )_{b^s_{2,2}} := \inf_{|\Lambda| \le n}\,  \Big\|
\sum_{(j,\lambda) \not\in \Lambda} \, a_{j, \lambda}\,  e_{j,\lambda}  
\, \Big\|_{b^s_{2,2}(\nabla)} \, .
\]
If $\Lambda$ contains the $n$ largest terms 
$2^{js}\, |\langle f, {\cal  E}^*\tilde{\psi}_{j,\lambda}\rangle | $ then
\[
\sigma_n (f,(\cF, \cG))_{B^s_2(L_2(\Omega))}
\le  c_1  \, \sigma_n \Big( (\langle f, {\cal  E}^*\tilde{\psi}_{j,\lambda}
  \rangle)_{(j,\lambda) \in {\nabla}}, \O \Big)_{b^s_{2,2}}
\]
follows. Next we shall use the following abbreviations:
let $F_1 = B^{s+t}_q (L_p (\Omega))$ and $F_2= b^{s+t}_{p,q}(\nabla)$.
Using Proposition \ref{discrete} with respect to $\nabla$ and a simple homogeneity argument we find
\[
\sup_{\| f \|_{F_1} \le 1} \, \sigma_n (f,(\cF, \cG))_{B^s_2(L_2(\Omega))}
\le c_2 \, \sup_{\| a \|_{F_2} \le 1} \,
 \sigma_n ( a, \O )_{b^s_{2,2}}
\le c_3 \,  n^{-t/d}\, , 
\]
since 
\[ 
\|\,  (\langle f, {\cal  E}^*\tilde{\psi}_{j,\lambda}
  \rangle)_{j, \lambda \in {\nabla}}\|_{b^{s+t}_{p,q}} \asymp
 \|\, f\, \|_{B^{s+t}_q(L_p(\Omega))}\, .
\]
This completes the proof of Theorem \ref{th5}.
{\eproof}

\begin{rem}
The advantage of our frame construction consists  in the fact
that it is universal for all bounded Lipschitz domains.
The disadvantage of our frame construction lies 
in the use of the operator $\ce^*$. 
This limits its value in case of concrete calculations.
There are other frame constructions in the 
literature. Let us mention here
the constructions given in {\rm \cite{CDD00}, \cite{T06b}} and 
{\rm \cite{DFR04}}. We add a few comments to these frames:
\begin{itemize}
\item
The frame pairs constructed 
in {\rm \cite{CDD00} }allow a discretization of Besov spaces
on domains $\Omega$ under certain restrictions, 
both with respect to the domains and with 
respect to the parameters of the Besov space.
In particular, only the case $1\le p\le \infty$, 
$0 < q \le \infty$ and $s>0$ is considered. 
With $(\cF,\cG)$ denoting the frame pairs constructed in the 
aforementioned paper we obtain
\[
\sup_{\| f \|_{F_1} \le 1} \, \sigma_n (f,(\cF, \cG))_{H^{-s}(\Omega)} \asymp n^{-t/d}
\]
where 
\[
F_1:=  B^{-s+t}_q (L_p (\Omega)) \, , \qquad t-s>0\, , \quad 1\le p,q\le \infty \, .
\]
Generalization to the case $0 <q,p <1$  have been given in {\rm \cite{DKW}}. 
\item
The frames constructed in {\rm \cite{T06b}} allow a discretization of Besov spaces
on Lipschitz domains $\Omega$ 
under the restrictions $0< p,q \le \infty$ and $s<0$.
The frame pairs  consist of either wavelets originating from a wavelet basis on $\Rd$
or dilated and shifted versions of the associated scaling function.
They all have the property that their support is contained in $\Omega$.
Furthermore, these dilated and shifted copies of the scaling functions show up only near the boundary. Inside a box contained in $\Omega$ and with some distance to the boundary
the frame pair reduces to a biorthogonal wavelet subsystem.
The same construction can be made to discretize the Besov spaces
$\widetilde{B}^s_q (L_p (\Omega))$ if $s> d \max(0, 1/p-1)$, see
the Appendix for a definition.
Hence, with $(\cF,\cG)$ denoting the frame pair of {\rm \cite{T06b}}  we obtain
\[
\sup_{\| f \|_{F_1} \le 1} \, 
\sigma_n (f,(\cF, \cG))_{H^{-s}(\Omega)} \asymp n^{-t/d}, 
\]
where 
\[
F_1:= \left\{\begin{array}{lll} B^{-s+t}_q (L_p (\Omega)) & 
\qquad & \mbox{if} \quad
t-s<0\\
\widetilde{B}^{-s+t}_q (L_p (\Omega)) && 
\mbox{if}\quad t-s>\, d \max(0,\frac 1p-1). 
             \end{array} \right.
\]
\item
The frame pairs constructed in {\rm \cite{DFR04}} allow a discretization of  $H^s(\Omega)$-spaces with $s>0$. This construction works for domains with piecewise analytic boundary and is based on an {\rm overlapping} partition of the
domain by means of sufficiently smooth parametric images of the unit cube. On the reference cube,  a tensor product  biorthogonal  wavelet basis
employing the  boundary adapted wavelets on the interval  from {\rm \cite{DS98}} is constructed.   Under certain conditions,
 the union of all the parametric images of these bases  gives rise to frame pair for $H^s(\Omega), ~s>0.$
 \item Of course, all the examples of biorthogonal wavelet bases on  polyhedral  domains also fit into our setting. One natural way as, e.g., outlined in {\rm \cite{CTU}} and
{\rm \cite{DS99}}, is to decompose the domain into a  {\rm disjoint}  union 
of parametric
images of reference cubes. Then  one constructs wavelet bases on the reference
cubes and glues everything together in a judicious fashion. However, due to the
glueing procedure, only Sobolev spaces $H^s$ with smoothness $s <3/2$ can be
characterized. This bottleneck can be circumvented by the approach in
{\rm \cite{DS99a}}.  
There,  a much more tricky domain decomposition method involving
certain projection and extension operators is used. By proceeding in 
this way, norm
equivalences for all spaces $B^t_q(L_p(\Omega))$ can be derived, at least for
the case $p>1$, see {\rm \cite[Theorem 3.4.3]{DS99a}}.  However, the 
authors also mention that their results can be generalized to the case $p<1$, 
see {\rm \cite[Remark 3.1.2]{DS99a}}.
\end{itemize}

\end{rem}


\subsection{Proof of Theorem \ref{basic2}}\label{basic3}


Periodic Besov spaces have analoguous properties than the Besov spaces defined on 
smooth domains or on $\Rd$. Our general reference for these classes is \cite{ScTr87}.
A definition of periodic Besov spaces is given in the Appendix.


\subsubsection{Widths of Periodic Besov Spaces}


As a preparation of the proof of Theorem \ref{basic2} we shall investigate 
the widths of embeddings of periodic Besov spaces, a topic which is also of self-contained interest.
In \cite{dns2} we reduced the corresponding problem for the nonperiodic 
Besov spaces on a Lipschitz domain to that one for the discrete Besov spaces.
It would be of interest to construct an isomorphism between 
these periodic spaces $B^{s}_{q} (L_p(\tor))$ and $b^s_{p,q}$ 
as well, see Subsection \ref{discr}.
Periodic wavelet constructions exist in the literature. However, up to our knowledge, those characterizations of periodic Besov spaces are established only with 
additional restrictions for the parameters. So we employ a different strategy here.

\begin{thm}   \label{basic4} 
Let $0 < p, q \le \infty$, $s \in \R $ 
and suppose that  
\[
t> \Big(\frac{1}{p} - \frac 12 \Big)_+
\] 
holds. Then 
there exists a constant $C^*$ such that for any $C \ge C^*$
we have 
\[
e_{n,C}^{\rm frame} (I, B^{s+t}_{q} 
(L_p(\tor)), B^s_{2} (L_2(\tor)))   \asymp
n^{-t}\,  . 
\]
\end{thm} 

\begin{proof}
{\em Step 1.} Preparations.
For the estimate from above we shall use a connection between periodic and weighted spaces.
Let $\varrho_\kappa (x):= (1+|x|^2)^{-\kappa/2}$, $x\in \R$, $\kappa>0$.
We define
\begin{equation}
B^s_q(L_p(\R, \varrho_\kappa)) := \Big\{f\in \cs'(\R):\quad 
f \, \varrho_\kappa \in B^s_q(L_p (\R))\Big\}\, ,
\end{equation}
endowed with the natural quasi-norm
\[
\| \, f \, |B^s_q(L_p (\R, \varrho_\kappa))\| := \| \, f \, \varrho_\kappa \, |
 B^s_q(L_p (\R))\| \, .
\]
Here $\cs'(\R)$ denotes the collection of the tempered distributions on $\R$.
As a combination of Franke's characterization of weighted spaces, see
Theorem 5.1.3 in \cite{ScTr87}, and a result of Triebel \cite{T80}
we find that $f \in B^{s}_{q} (L_p(\tor))$ if and only if
$f$ is a $2\pi$-periodic distribution in $\cs'(\R)$ which belongs to 
$B^s_q(L_p (\R, \varrho_\kappa))$ with $\kappa> (1/p)$.
Moreover, there exist positive constants $c_1,c_2$ such that
\[
c_1\, \| \, f \, |B^s_q(L_p (\R, \varrho_\kappa))\| \le 
\| \, f \, |B^s_q(L_p (\tor)\|\le c_2\, 
\| \, f \, |B^s_q(L_p (\R, \varrho_\kappa))\|
\]
holds for all such $f$.
\\
{\em Step 2.} Let $\psi \in C_0^\infty (\R)$ be a smooth cut-off function
such that  $\psi (x)=1$ if $|x|\le \pi$ and $\psi (x)=0$ if $|x|\ge 2\pi$. We shall study the mapping $T: \, f \mapsto \psi \cdot f$. Let $J=[-3\pi,3\pi]$.
Obviously
\begin{eqnarray*}
\| \, f \, \psi \, |B^s_q(L_p (J))\|& \le &  \| \, f \, \psi \, |B^s_q(L_p (\R))\|
= \| \, f \, \psi \, \varrho_\kappa \, (1/\varrho_\kappa) \,
\, \psi (\cdot/2)\,  |B^s_q(L_p (\R))\|
\\
& \le & c_3 \, \|\, (1/\varrho_\kappa) \,
\, \psi (\cdot/2)\,  |C^\mu (\R)\|\, 
\| \, f \, \psi \, \varrho_\kappa \,  |B^s_q(L_p (\R))\|
\\
& \le & c_4 \, 
\| \, f \, \psi \,   |B^s_q(L_p (\R,\varrho_\kappa))\|\, .
\end{eqnarray*}
where $\mu$ has to be chosen sufficiently large, cf. e.g. \cite[2.8]{T83} or
\cite[4.7]{RS96}.
Since $\psi$ is a pointwise multiplier for these weighted Besov spaces as well
we end up with $T \in \cl (B^{s}_{q} (L_p(\tor)), B^{s}_{q} (L_p(J)))$.
Moreover, $T$ is a bijection onto a closed subspace of 
$B^{s}_{q} (L_p(J))$, denoted by $T^{s}_{q} (L_p(J))$,
simultenuously for all parameters.
Now we consider the commutative diagram:
\[
\begin{CD}
B^{s+t}_{q}(L_p(\tor))@> I_1 >> B^{s+t}_{2} (L_2(\tor))\\
@VTVV @AA T^{-1} A\\
T^{s+t}_{q} (L_p(J)) @> I_2 >> T^{s}_{2} (L_2(J)) \, .
\end{CD}
\]
Lemma \ref{zusatz} yields 
\[
e_{n,\wt C}^{\rm frame} (I_1, B^{s+t}_{q} 
(L_p(\tor)), B^s_{2} (L_2(\tor)))   \le 
 \| T\| \, \| T^{-1} \|  \, 
e^{\rm frame}_{n,C} (I_2,T^{s+t}_{q} (L_p(J)),T^{s}_{2} (L_2(J)))
\]
with $\wt C = C \, \Vert T^{-1} \Vert \, \Vert T \Vert$.
Now we employ Lemma \ref{zusatz2} and obtain
\[
e^{\rm frame}_{n,C} (I_2,T^{s+t}_{q} (L_p(J)),T^{s}_{2} (L_2(J)))
\le  e^{\rm frame}_{n,C} (I_2,T^{s+t}_{q} (L_p(J)),B^{s}_{2} (L_2(J)))\, .
\]
This, together with a monotonicity arguments leads to
\[
e_{n,\wt C}^{\rm frame} (I_1, B^{s+t}_{q} 
(L_p(\tor)), B^s_{2} (L_2(\tor)))   \le 
 \| T\| \, \| T^{-1} \|  \, 
e^{\rm frame}_{n,C} (I_2,B^{s+t}_{q} (L_p(J)),B^{s}_{2} (L_2(J)))\, .
\]
The estimate from above is finished by using Theorem \ref{th5} with $\Omega=J$
and $d=1$.\\
{\em Step 3.} 
Let $J=(-1/2,1/2)$. Then there exists a linear extension operator 
$\ce : \, B^s_q (L_p (J)) \to B^s_q(L_p (\R))$, see \cite{Ry}.
Let $\psi$ be as above.
We define
\[
Tf (x) := \left\{
\begin{array}{lll}
\ce f (x) \, \psi (6x) && \mbox{if} \quad -\pi \le x \le \pi\, , \\
2\pi\mbox{-periodic extension} && \mbox{otherwise}\, .
\end{array}
\right. 
\]
We claim that $T \in \cl (B^s_q (L_p(J)), B^s_q (L_p (\tor)))$ for all parameter constellations.
To see that we first construct an appropriate 
decomoposition of unity. We put
\[
\varphi (x):= \frac{\psi (x)}{\sum_{k=-\infty}^\infty \psi (x-2\pi k)}\, , 
\quad x \in \R\, .
\]
It follows that 
\[
1= \sum_{m=-\infty}^\infty \varphi (x-2\pi m) \qquad \mbox{for all}\quad x \in \R
\]
and $\supp \varphi \subset \{x \in \R: \: \psi (x/2)=1\}$.
Hence, with $t= \min(1,p,q)$ and $\kappa > 1/t \ge 1/p$, we obtain
\begin{eqnarray*}
&& \hspace*{-0.8cm}
\| \, Tf \, |B^s_q (L_p (\tor))\|^t  \le  c_2^t \, \| \, (Tf) \, \varrho_\kappa \, |B^s_q (L_p (\R))\|^t
\\
& = & c_2^t \, \| \, \sum_{m=-\infty}^\infty  \varphi (\cdot -2\pi m)\, (Tf) \, \varrho_\kappa \, |B^s_q (L_p (\R))\|^t  \\
& \le & c_2^t \,  \sum_{m=-\infty}^\infty  \| \, \varphi (\cdot -2\pi m)\, (Tf) \, \varrho_\kappa \, |B^s_q (L_p (\R))\|^t  \\
& = & c_2^t \,  \sum_{m=-\infty}^\infty  \| \, \varphi (\cdot -2\pi m)\, 
\psi \Big(\frac{\cdot -2\pi m}{2}\Big)\,
(Tf) \, \varrho_\kappa \, |B^s_q (L_p (\R))\|^t  \\
& \le & c_3 \,  \sum_{m=-\infty}^\infty  \| \, \varphi (\cdot -2\pi m)\,  \varrho_\kappa \, |C^\mu (\R)\|^t \, 
 \| \, \psi \Big(\frac{\cdot -2\pi m}{2}\Big)\,  (Tf)  \, |B^s_q (L_p (\R))\|^t  \, ,
\end{eqnarray*}
where we used again assertions on pointwise multipliers, see, e.g., 
\cite[2.8]{T83} or \cite[4.7]{RS96}.
The shift-invariance of $\| \, \cdot \, |B^s_q (L_p (\R))\|$ and the periodicity of $Tf$
imply
\[
\| \, \psi \Big(\frac{\cdot -2\pi m}{2}\Big)\,  (Tf)  \, |B^s_q (L_p (\R))\|
=
\| \, \psi (\, \cdot \, /2)\,  (Tf)  \, |B^s_q (L_p (\R))\|
\]
for all $m \in \Z$.
Furthermore, elementary calculations yield
\[
\| \, \varphi (\cdot -2\pi m)\,  \varrho_\kappa \, |C^\mu (\R)\| 
\le c_4\, \varrho_\kappa (2\pi m) 
\]
with $c_4$ independent of $m$. Altogether this proves
\begin{eqnarray*}
\| \, Tf \, |B^s_q (L_p (\tor))\|
& \le &  c_5 \, 
\| \, \psi (\, \cdot \, /2)\,  (Tf)  \, |B^s_q (L_p (\R))\|
\Big( \sum_{m=-\infty}^\infty \varrho (2\pi m)^t  \Big)^{1/t}
\\
& \le &  c_6 \, 
\| \, \psi (\, \cdot \, /2)\,  (Tf)  \, |B^s_q (L_p (\R))\|\, .
\end{eqnarray*}
Taking into account the identity
\[
\psi (x/2)\,  Tf (x) = \psi (x/2)\, \Big(\sum_{m=-2}^2
\ce f (x-2\pi m) \, \psi (6(x-2\pi m)) \Big) 
\]
we have  
\begin{eqnarray*} 
\| \, \psi (\, \cdot \, /2)\,  (Tf)  \, |B^s_q (L_p (\R))\| & \le & c_7 \sum_{m=-2}^2
\| \, \psi (\, \cdot \, /2)\,  \ce f (x-2\pi m) \, \psi (6(x-2\pi m)) \, |B^s_q (L_p (\R))\| \\
& \le & c_8 \sum_{m=-2}^2
\| \,  \ce f (x-2\pi m) \, \psi (6(x-2\pi m)) \, |B^s_q (L_p (\R))\| \\
& \le & c_9 
\| \,  \ce f  \, \psi (6(\, \cdot\, )) \, |B^s_q (L_p (\R))\| \\ 
& \le & c_{10} 
\| \,  \ce f \, |B^s_q (L_p (\R))\| \\ 
& \le & c_{10} \, \| \ce \|\, \| \,  f \, |B^s_q (L_p (J))\| \, ,
\end{eqnarray*}
which proves the claim.
Moreover, $T$ is a bijection onto a closed subspace of $B^s_q (L_p (\tor))$.
This subspace will be denoted by $T^s_q (L_p (\tor))$.
Now we can argue as in Step 2. The commutative diagram
\[
\begin{CD}
B^{s+t}_{q}(L_p(J))@> I_1 >> B^{s+t}_{2} (L_2(J))\\
@VTVV @AA T^{-1} A\\
T^{s+t}_{q} (L_p(\tor)) @> I_2 >> T^{s}_{2} (L_2(\tor)) 
\end{CD}
\]
implies 
\[
e_{n,\wt C}^{\rm frame} (I_1, B^{s+t}_{q} 
(L_p(J)), B^s_{2} (L_2(J)))   \le 
 \| T\| \, \| T^{-1} \|  \, 
e^{\rm frame}_{n,C} (I_2,B^{s+t}_{q} (L_p(\tor)),B^{s}_{2} (L_2(\tor)))\, .
\]
with $\wt C = C \, \Vert T^{-1} \Vert \, \Vert T \Vert$.
The estimate from below is finished by using Theorem \ref{th5} with $\Omega=J$
and $d=1$.
\end{proof}

Now we consider some subspaces of $B^s_q (L_p (\tor))$.
Let
\begin{equation}\label{eq-61}
Z^s_q (L_p (\tor)) := \Big\{f \in B^s_q (L_p (\tor)): \quad
\langle f,1 \rangle_\tor =0 \Big\}\, .
\end{equation}
Observe that the function $g(x)=1$ belongs to $D(\tor)$, the collection of all 
complex-valued, $2\pi$-periodic and infinitely differentiable function.
Since
\[
D(\tor) \hookrightarrow B^s_q (L_p (\tor)) \hookrightarrow D' (\tor)
\]
the scalar product $\langle f,1 \rangle_\tor $ is well-defined for all $f \in
B^s_q (L_p (\tor)) $, cf. \cite[3.5.1]{ScTr87}.

\begin{cor}   \label{basic5} 
Let $0 < p, q \le \infty$, $s \in \R $ 
and suppose that  
\[
t> \Big(\frac{1}{p} - \frac 12 \Big)_+
\] 
holds. Then 
there exists a constant $C^*$ such that for any $C \ge C^*$
we have 
\[
e_{n,C}^{\rm frame} (I, Z^{s+t}_{q} 
(L_p(\tor)), Z^s_{2} (L_2(\tor)))   \asymp
n^{-t}\,  . 
\]
\end{cor}

\begin{proof}
The upper estimate  can be established as above.
For the estimate from below 
we start with $f \in B^s_q (L_p (J))$ and $J=[-1/2,-1/4]$.
The operator $T$ has to be replaced by 
\[
\wt T f (x):= \left\{
\begin{array}{lll}
\ce f (x) \, \psi (14(x+1/2)) -
\ce f (-x) \, \psi (14(-x+1/2)) 
 && \mbox{if} \quad -\pi \le x \le \pi\, , \\
2\pi\mbox{-periodic extension} && \mbox{otherwise}\, .
\end{array}
\right.
\]
Hence $\langle \wt T f, 1 \rangle_\tor = 0 $ which is clear for $f\in D(\tor)$.
Since $D(\tor)$ is dense in $D'(\tor)$ it follows in general.
\end{proof}


\subsubsection{Besov Spaces on the Unit Circle}\label{circle}


There is a simple transformation of the interval $[0,2\pi)$ onto the unit circle
given by
\[
t \mapsto (\cos t, \sin t)\, , \qquad 0 \le t < 2\pi\, .
\]
For a given distribution $f \in D'(\Gamma)$  
we define
\begin{equation}
h (t):= f (\cos t, \sin t)\, , \qquad t \in \R\, .
\end{equation}
Observe that $\varphi \in D(\Gamma)$ implies
$\varphi (\cos t, \sin t ) \in D(\tor)$. Hence, if $f \in D'(\Gamma)$
then $h \in D'(\tor)$.

\begin{definition}
Let $s \in \R$ and $0 < p,q \le \infty$. Then 
$B^s_q(L_p (\Gamma))$ is the collection of all 
distributions $f \in D'(\Gamma)$ such that
the corresponding distribution $h$ is contained in $B^s_q(L_p (\tor))$.
We put
\[
\| \, f \, | B^s_q(L_p (\Gamma))\|:= 
\| \, h \, | B^s_q(L_p (\tor))\|.
\]
\end{definition}

\begin{lemma}
In the sense of equivalent norms we have
$H^{1/2}(\Gamma) = B^{1/2}_2 (L_2 (\Gamma))$ as well as
$H^{-1/2}(\Gamma) = B^{-1/2}_2 (L_2 (\Gamma))$.
\end{lemma}

\begin{proof}
It holds
\[
 B^{1/2}_2 (L_2 (\tor)) = \Big\{
h \in L_2 (\tor): \quad 
\int_0^{2\pi} \int_0^{2\pi} \frac{|h(x) - h(y)|^2}{|x-y|^2}\, dx\, dy < \infty
\Big\}\, , 
\]
see e.g. \cite[3.5.4]{ScTr87}.
Furthermore, the norms $\| \, h \, | B^s_q(L_p (\tor))\|$ and 
\[
\| \, h \, | L_2 (\tor)\| + \Big( 
\int_0^{2\pi} \int_0^{2\pi} \frac{|h(x) - h(y)|^2}{|x-y|^2}\, dx\, dy\Big)^{1/2} 
\]
are equivalent. Now it remains to observe that
\begin{eqnarray*}
\| \, f \, | L_2 (\Gamma)\| &  + & \Big( 
\int_\Gamma \int_\Gamma \frac{|f(x) - f(y)|^2}{|x-y|^2}\, d\Gamma_x\, d\Gamma_y\Big)^{1/2} 
\\
& \asymp &
\| \, h \, | L_2 (\tor)\| + \Big( 
\int_0^{2\pi} \int_0^{2\pi} \frac{|h(x) - h(y)|^2}{|x-y|^2}\, dx\, dy\Big)^{1/2}
\end{eqnarray*}
since there exist positive constants $c_1, c_2$ such that
\[
c_1 \, |x-y|^2 \le (\cos x - \cos y)^2 + (\sin x - \sin y)^2 \le c_2 \, |x-y|^2
\]
for all $x,y \in [0,2\pi]$.
This proves $H^{1/2}(\Gamma) = B^{1/2}_2 (L_2 (\Gamma))$ in the sense of 
equivalent norms.
The second assertion follows from
$(H^{1/2}(\Gamma))' = H^{-1/2}(\Gamma)$
(just by definition) and the duality relation
$(B^{1/2}_2 (L_2 (\tor)))' =  B^{-1/2}_2 (L_2 (\tor))$, see \cite[3.5.6]{ScTr87}.
\end{proof}


\subsubsection{Proof of Theorem \ref{basic2}}\label{basic6}


We consider the commutative diagram
\[
\begin{CD}
Y^{t+1/2}_{q}(L_p(\Gamma))@> I_1 >> H^{1/2} (\Gamma)\\
@VTVV @AA T^{-1} A\\
Z^{t+1/2}_{q} (L_p(\tor)) @> I_2 >> Z^{1/2}_{2} (L_2(\tor)) 
\end{CD}
\]
Here the operator $T$ is chosen to be the mapping $f \mapsto h$.
Since $T$ is a bijection considered as a mapping defined on $D'(\Gamma)$
with values in $D'(\tor)$ we obtain that $T$ is an isomorphism 
belonging to $\cl (B^{t+1/2}_{q}(L_p(\Gamma)), B^{t+1/2}_{q}(L_p(\tor)))$.
Consequently, $T:\,  Y^{t+1/2}_{q}(L_p(\Gamma)) \to Z^{t+1/2}_{q}(L_p(\tor))$
is an isomorphism as well.
Lemma \ref{zusatz} yields
\begin{equation}\label{eq-62}
e^{\rm frame}_{n, \wt C} (I_1,Y^{t+1/2}_{q}(L_p(\Gamma)),H^{1/2} (\Gamma)) \le  
\| T\| \, \| T^{-1} \|  \, 
e^{\rm frame}_{n,C} (I_2,Z^{t+1/2}_{q}(L_p(\tor)),Z^{1/2}_{2}(L_2(\tor)))
\end{equation} 
with $\wt C = C \, \Vert T^{-1} \Vert \, \Vert T \Vert$.
As a consequence of the commutative diagram
\[
\begin{CD}
Z^{t+1/2}_{q}(L_p(\tor))@> I_1 >> Z^{1/2}_2 (L_2(\tor))\\
@VT^{-1}VV @AA T A\\
Y^{t+1/2}_{q} (L_p(\Gamma)) @> I_2 >> H^{1/2} (\Gamma) 
\end{CD}
\]
Lemma \ref{zusatz}, and inequality (\ref{eq-62}) we conclude
\[
e^{\rm frame}_{n, \wt C} (I_1,Y^{t+1/2}_{q}(L_p(\Gamma)),H^{1/2} (\Gamma)) \asymp
e^{\rm frame}_{n,C} (I_2,Z^{t+1/2}_{q}(L_p(\tor)),Z^{1/2}_{2}(L_2(\tor)))\, .
\]
>From Corollary \ref{basic5} we derive
\[
e^{\rm frame}_{n, C} (I_1,Y^{t+1/2}_{q}(L_p(\Gamma)),H^{1/2} (\Gamma)) \asymp n^{-t}
\]
for $C$ sufficiently large.
Now the assertion follows from 
the commutative diagram (\ref{final}) and Lemma \ref{basic}.


\section{Appendix -- Besov Spaces}\label{Appendix}


Here we collect some properties of Besov spaces 
which have been used in the text before. 
For general information on Besov spaces we refer 
to the monographs \cite{Me,Ni,Pe,RS96,T83,T92,T06a}. A collection of results 
for Besov as well as Sobolev spaces on domains can be found  in \cite{dns2}.
There detailed references are given.\\
In most of the references given above Besov as well as Sobolev spaces are 
treated as classes of complex-valued functions (distributions).
In the framework of information based complexity it is common to deal with real-valued functions (distributions), cf. e.g. (\ref{manifold}).
Here we make use of the following point of view:
all spaces in the Appendix are spaces of complex-valued distributions.
Then, finally we consider the restrictions to the real-valued subspaces.


\subsection{Wavelet Characterizations}
\label{app3}


For the  construction of  biorthogonal 
wavelet  bases  as considered below
we refer to the recent monograph of 
Cohen \cite[Chapt.~2]{C03}. 
Let $\varphi$ be a compactly supported 
scaling function of sufficiently high regularity and let 
$\psi_i$, $i=1, \ldots \, 2^d-1$ be corresponding wavelets.
More exactly, we suppose for some $N>0$ and $r\in \Nb$
\begin{eqnarray*}
&& \supp \, \varphi\, , \: \supp \, \psi_i 
\quad \subset \quad  [-N,N]^d\, , 
\qquad i=1, \ldots , 2^d-1 \, , \\
&& \varphi,  \psi_i \in C^r(\Rd)\, , 
\qquad i=1, \ldots , 2^d-1 \, ,   \\
&&
\int x^\alpha \, \psi_i (x)\, dx =0 \qquad \mbox{for all} \quad
|\alpha|\le r\, , \qquad i=1, \ldots , 2^d-1 \, ,  
\end{eqnarray*}
and
\[
\varphi (x - k), \, 2^{jd/2}\, \psi_i (2^jx-k)\, , 
\qquad j\in \Nb_0\, ,
\quad k \in \Zd \, , 
\]
is a Riesz basis in $L_2(\Rd)$.
We shall use the standard abbreviations
\[
\psi_{i,j,k} (x) = 2^{jd/2}\, \psi_i (2^jx-k) \qquad \mbox{and} \qquad
\varphi_k (x) = \varphi (x-k) \, .
\]
Further, 
the dual Riesz basis should fulfill  the same requirements, i.e., 
there exist functions $\widetilde{\varphi} $ and 
$\widetilde{\psi}_i$, $ i=1, \ldots , 2^d-1$,
such that 
\begin{eqnarray*}
\langle \widetilde{\varphi}_k, \psi_{i,j,k} \rangle   
& = & \langle \widetilde{\psi}_{i,j,k}, 
\varphi_k \rangle  = 0\, , \\
\langle \widetilde{\varphi}_k, \varphi_{\ell} \rangle  
& = & \delta_{k,\ell} \qquad
(\mbox{Kronecker symbol})\, , \\
\langle \widetilde{\psi}_{i,j,k}, 
\psi_{u,v,\ell} \rangle   & = & \delta_{i,u}\,
 \delta_{j,v}\, \delta_{k,\ell}\, , \\
\supp \, \widetilde{\varphi} \, , &&  
\supp \, \widetilde{\psi}_i \quad 
\subset \quad  [-N,N]^d\, , 
\qquad i=1, \ldots , 2^d-1 \, , \\
\widetilde{\varphi},   \widetilde{\psi}_i & \in 
&  C^r(\Rd)\, , \qquad i=1, \ldots , 2^d-1 \, ,   \\
\int x^\alpha \, 
\widetilde{\psi}_i (x)\, dx & = & 0 \qquad \mbox{for all} \quad
|\alpha|\le r\, , \qquad i=1, \ldots , 2^d-1 \, .  
\end{eqnarray*}  
For $f \in \S'(\Rd)$ we put
\begin{equation}\label{eq100}
\langle f,\psi_{i,j,k} \rangle  
= f(\overline{\psi_{i,j,k}}) \qquad \mbox{and}\qquad
 \langle f,\varphi_k \rangle  = f(\overline{\varphi_k}) \, ,
\end{equation}
whenever this makes sense.

\begin{proposition}\label{wavelets}
Let $s \in \R$ and $0 < p,q \le \infty$. 
Suppose 
\begin{equation}\label{eq101}
r > \max \Big (s, d\, \max (0, \frac 1p -1) - s\Big)\, .
\end{equation}
Then $B^s_q(L_p(\Rd))$ is the collection of all 
tempered distributions $f$ such that
$f$ is representable as
\[
f = \sum_{k \in \Zd} a_{k}\, \varphi_k
 + \sum_{i=1}^{2^d-1} \, \sum_{j=0}^\infty 
\sum_{k \in \Zd} a_{i,j,k}\, \psi_{i,j,k} 
\qquad (\mbox{convergence in} \quad \S')
\]
with
\[
\| \, f\, |B^s_q(L_p(\Rd))\|^* 
:= \Big(\sum_{k \in \Zd} |a_{k}|^p\Big)^{1/p} +
 \bigg(\sum_{i=1}^{2^d-1} \, \sum_{j=0}^\infty 
2^{j(s + d(\frac{1}{2}-\frac{1}{p}))q} 
\Big(\sum_{k \in \Zd} |a_{i,j,k}|^p\Big)^{q/p}\bigg)^{1/q} < \infty\, ,
\]
if $q < \infty$ and
\[
\| \, f\, |B^s_\infty (L_p(\Rd))\|^* := 
\Big(\sum_{k \in \Zd} |a_{k}|^p\Big)^{1/p} +
  \sup_{i=1, \ldots \, , 2^d-1} \, \sup_{j=0, \ldots \, } \, 
2^{j(s + d(\frac{1}{2}-\frac{1}{p}))} 
\Big(\sum_{k \in \Zd} |a_{i,j,k}|^p\Big)^{1/p} < \infty\, .
\]
The representation is unique and
\[
a_{i,j,k} =  \langle f,
\widetilde{\psi}_{i,j,k}\rangle  
\qquad \mbox{and} \qquad 
a_k = \langle f,\widetilde{\varphi}_{k} \rangle  
\] 
hold. Further $I : f \mapsto \{\langle f,
\widetilde{\varphi}_{k} \rangle ,\,  \langle f,
\widetilde{\psi}_{i,j,k}\rangle  \} $
is an isomorphic map of $B^s_q(L_p (\Rd))$ onto 
the sequence space (equipped with the quasi-norm 
$\| \, \cdot \, |B^s_q(L_p(\Rd))\|^*$), i.e. 
$\| \, \cdot \, |B^s_q(L_p(\Rd))\|^*$ may serve as an 
equivalent quasi-norm on $ B^s_q(L_p (\Rd))$.
\end{proposition}

A proof of Proposition \ref{wavelets} has been given 
in {\rm \cite{T06b}},  see also \cite{Ky03}
for a homogeneous version.
A different proof, but restricted to 
$s> d(\frac 1p - 1)_+$, is given in {\rm \cite[Thm.~3.7.7]{C03}.}
However, there are many forerunners with 
some restrictions on $s,p$ and $q$.


\subsection{Besov Spaces on  Domains}\label{app4}


Let $\Omega \subset \Rd$ be an  bounded 
open nonempty set. Then we define
$B^s_q (L_p(\Omega))$ to be the collection of 
all distributions $f \in \D'(\Omega)$ such that there
exists a tempered distribution $g \in B^s_q (L_p (\Rd))$ satisfying 
\[
f (\varphi) = g(\varphi) \qquad \mbox{for all} 
\quad \varphi \in \D(\Omega) \, ,
\]
i.e. $g|_\Omega = f$ in $\D'(\Omega)$.
We put
\[
\| \, f\, |B^s_q (L_p (\Omega))\|:= 
\inf \, \| \, g \, |B^s_q (L_p (\Rd))\|\, , 
\] 
where the infimum is taken with respect 
to all distributions $g$ as above.

%
%

\subsection{Sobolev Spaces on Domains}\label{app6}


Let $\Omega $ be a bounded Lipschitz domain.
Let $m \in \Nb$.
As usual $H^m (\Omega)$ denotes the 
collection of all functions $f$ such that the 
distributional derivatives $D^\alpha f$ of order  $|\alpha|\le m$
belong to $L_2 (\Omega)$. The norm is defined as
\[
\| \, f \, |H^m (\Omega)\| 
:= \sum_{|\alpha|\le m} \| \, D^\alpha f\, |L_2(\Omega)\|\, .
\]
It is well-known that $H^m (\Rd) = B^m_2(L_2(\Rd))$ 
in the sense of equivalent norms,
cf. e.g. \cite{T83}.
As a consequence of the existence of a 
bounded linear extension operator 
for Sobolev spaces on bounded Lipschitz 
domains, cf. \cite[p.~181]{St}, it follows
\[
H^m (\Omega) = B^m_2(L_2(\Omega)) \qquad \mbox{(equivalent norms)} \, ,
\]
for such domains. For fractional $s>0$ we introduce the classes by 
complex interpolation. 
Let $0 < s < m$, $s \not\in \Nb$. 
Then, following \cite[9.1]{LM}, we define
\[
H^s (\Omega) := \Big[H^m (\Omega), L_2(\Omega) \Big]_\Theta\, , 
\qquad \Theta = 1-\frac sm\, .
\]
This definition does not depend on $m$ in 
the sense of equivalent norms, cf.
\cite{T02}. The outcome $H^s (\Omega)$ 

coincides with $ B^s_2(L_2(\Omega))$, 
cf. \cite{dns2} for further details.


\subsection{Spaces on Domains and Boundary Conditions}
\label{app7}


We concentrate on homogeneous boundary conditions.
Here it makes sense  to introduce two further 
scales of function spaces (distribution spaces).

\begin{definition}
Let $\Omega \subset \Rd$ be an open nontrivial set. 
Let $s \in \R$ and $0 < p,q \le \infty$. 
\\
{\rm (i)}
Then 
$\mathring{B}^s_q(L_p (\Omega))$ denotes the 
closure of ${\mathcal D}(\Omega)$
in $ B^s_q(L_p (\Omega))$, 
equipped with the quasi-norm of  $ B^s_q(L_p (\Omega))$.
\\
{\rm (ii)} Let $s\ge 0$.
Then 
${H}^s_0 (\Omega)$ denotes the closure of ${\mathcal D}(\Omega)$
in $ H^s(\Omega)$, equipped with the norm of  $ H^s(\Omega)$.
\\
{\rm (iii)} By 
$\widetilde{B}^s_{q}(L_p (\Omega))$ we denote the collection of all
$f \in  {\mathcal D}' (\Omega)$ such that 
there is a $g \in B^s_q(L_p (\Rd))$ with
\begin{equation}\label{rand}
g_{|_{\Omega}} = f \qquad \mbox{and}
\qquad \supp g \subset \overline{\Omega}\, ,
\end{equation}
equipped with the quasi-norm 
\[
\| \, f \, |   
\widetilde{B}^s_q(L_p (\Omega))\| 
= \inf \| \, g \, | B^s_q(L_p (\Rd))\|\, ,
\]
where the infimum is taken over all such distributions 
$g$ as in {\rm (\ref{rand}).}
\end{definition}

\begin{rem}\label{noch eins}
For a bounded Lipschitz domain it holds 
$\mathring{B}^s_q(L_p (\Omega)) 
=\widetilde{B}^s_q(L_p (\Omega)) = B^s_q(L_p(\Omega))$
if 
\[
 0 <p,q < \infty \, , \quad \max \Big(\frac 1p - 1, \, 
d \, \Big(\frac 1p - 1\Big)\Big) < s < \frac 1p\, , 
\]
cf. {\rm  \cite[Cor.~1.4.4.5]{Gr3}} and {\rm \cite{T02}}. Hence,
\[
H^s_0 (\Omega) = \mathring{B}^s_2(L_2(\Omega)) 
= \widetilde{B}^s_2(L_2 (\Omega)) = {B}^s_2(L_2 (\Omega)) =
H^s (\Omega)
\]
if $0 \le s <1/2$.
\end{rem}


\subsection{Sobolev Spaces with Negative Smoothness} \label{negative}


In what follows duality has to be understood in the framework of the dual pairing
$({\mathcal D}(\Omega), {\mathcal D}'(\Omega))$.

\begin{definition}
Let $\Omega \subset \Rd$ be a bounded Lipschitz domain.
For $s>0$ we define
\[
H^{-s} (\Omega) :=\left\{ \begin{array}{lll}
\Big({H}^s_0 (\Omega)\Big)' & \qquad & \mbox{if}
\quad  s-\frac 12 \neq \mbox{integer}\, , \\
&& \\
 \Big(\widetilde{B}^s_2 (L_2(\Omega))\Big)' && \mbox{otherwise}\, .
\end{array}\right.
\]
\end{definition}

\begin{rem}\label{tilde}
If $\Omega \subset \Rd$ is a bounded Lipschitz domain then 
\[
{H}^s_0 (\Omega) = \widetilde{B}^s_2 (L_2(\Omega))\, , 
\qquad s>0\, , \quad s-\frac 12 \neq \mbox{integer}\, ,
\]
holds. Furthermore
\begin{equation}
H^{-s} (\Omega) = B^{-s}_2 (L_2 (\Omega))\, , \qquad s>0\, , 
\end{equation}
to be understood in the sense of equivalent norms. 
Again we refer to {\rm \cite{dns2}} for detailed references.
\end{rem}


\subsection{Besov Spaces on the Torus} \label{periodic}


Here our general reference is \cite[Chapt.~3]{ScTr87}. Since we are using also spaces with negative smoothness 
$s<0$ and/or $p,q<1$  we shall give a 
definition, which relies on  Fourier analysis.

Let $D(T) $ denote the collection of all complex-valued infinitely differentiable functions on $\tor$ (i.e. $2\pi$-periodic).
By $D'(T)$ we denote its dual. 
Any $f \in D'(T)$ can be identified with its Fourier series
$\sum_{k=-\infty}^\infty c_k (f)\, e^{ikx}$ where
$c_k (f)= (2\pi)^{-1}\, f(e^{-ikx})$. 
\\
Next we need a smooth dyadic decompositions of unity.
Let $\phi \in C^\infty_0 (\R)$ be a function 
such that $\phi (x)=1$ if $|x|\le 1$ and 
$\phi (x)=0$ if $|x|\ge 2$. Then we put
\begin{equation}
\phi_0 (x):= \phi(x), \qquad \phi_j (x) 
:= \phi (2^{-j}x) - \phi (2^{-j+1}x)\, , \quad j \in \Nb\, . 
\end{equation}
It follows
\[
\sum_{j=0}^\infty  \phi_j (x) = 1\, , \qquad x \in \R\, , 
\]
and \[
\supp \, \phi_j \subset \Big\{x \in \Rd: \quad 2^{j-2} 
\le |x| \le 2^{j+1}\Big\}\, , 
\qquad j=1,2, \ldots \, .
\] 
By means of these functions we define the Besov classes.

\begin{definition}
Let $s \in \R$ and $0 < p,q \le \infty$. Then 
$B^s_q(L_p (\tor))$ is the collection of all 
periodic tempered distributions $f$ such that
\[
\| \, f \, | B^s_q(L_p (\tor))\|= \bigg( \sum_{j=0}^\infty 
2^{sjq}\, 
\|  \sum_{k=-\infty}^\infty  \phi_j (k) \, c_k (f) \, e^{ikx}\, |L_p (\tor)\|^q\bigg)^{1/q} <\infty
\]
if $q< \infty$ and
\[
\| \, f | B^s_\infty(L_p (\tor))\|=  \sup_{j=0,1, \ldots}\, 2^{sj}\,  \|  
\sum_{k=-\infty}^\infty  \phi_j (k) \, c_k (f) \, e^{ikx}\, |L_p (\tor)\|
 <\infty
\]
if $q= \infty$. 
\end{definition}

\begin{rem}\label{confusion}
\begin{itemize}
\item[i)]
These classes are quasi-Banach spaces. 
They do not depend on the chosen function $\phi$
(up to equivalent quasi-norms). 
\item[(ii)]
There is a number of different characterizations of periodic Besov spaces, cf. e.g.
{\rm \cite[Chapt.~3]{ScTr87}.} 
In particular we wish to refer to the characterization by differences {\rm \cite[3.5.4]{ScTr87}.}
\end{itemize}
\end{rem}

\medskip
\noindent
{\bf Acknowledgment:} \
We thank Hans Georg  Feichtinger,  Massimo Fornasier and Hans Triebel for valuable 
remarks and comments that improved our paper.

\bigskip
\vbox{\noindent Stephan Dahlke\\
Philipps-Universit\"at Marburg\\
FB12 Mathematik und Informatik\\
Hans-Meerwein Stra\ss e\\
Lahnberge\\
35032 Marburg, Germany\\
e-mail: {\tt dahlke@mathematik.uni-marburg.de}\\
WWW: {\tt http://www.mathematik.uni-marburg.de/$\sim$dahlke/}\\}

\bigskip
\smallskip
\vbox{\noindent Erich Novak, Winfried Sickel\\
Friedrich-Schiller-Universit\"at Jena\\ 
Mathematisches Institut\\
Ernst-Abbe-Platz 2\\ 
07743 Jena, Germany\\
e-mail: {\tt \{novak, sickel\}@math.uni-jena.de}\\
WWW: {\tt http://www.minet.uni-jena.de/$\sim$\{novak,sickel\}}}

\bigskip

\end{document}